\documentclass[12pt,draftcls,onecolumn]{IEEEtran}

\IEEEoverridecommandlockouts                              

\makeatletter
\let\NAT@parse\undefined
\makeatother

\usepackage{graphicx}
\usepackage{bm}
\usepackage{epsfig}
\usepackage{amsthm}
\usepackage{amssymb,amsfonts,amsmath}
\usepackage{subfigure}
\usepackage{setspace}
\usepackage{enumerate}
\usepackage[numbers]{natbib}
\usepackage{color}
\usepackage{dsfont}
\usepackage{verbatim}

\newtheorem{theorem}{Theorem}
\newtheorem{lemma}{Lemma}
\newtheorem{remark}{Remark}
\newtheorem{definition}{Definition}
\newtheorem{proposition}{Proposition}
\newtheorem{example}{Example}

\newtheorem{corollary}{Corollary}

\def \w {\omega}
\def \t {\theta}
\def \a {\alpha}
\def \b {\beta}

\def \g {\gamma}
\def \l {\lambda}
\def \e {\varepsilon}
\def \< {\langle}
\def \> {\rangle}
\def \d {\delta}
\def \O {\Omega}
\def \ad{{\rm ad}}

\def \F{\mathcal{F}}
\def \B{\mathcal{B}}
\def \~{\tilde}

\usepackage{graphicx}
\usepackage{bm}
\usepackage{epsfig}
\usepackage{amssymb,amsfonts,amsmath}
\usepackage{subfigure}
\usepackage{setspace}
\usepackage{enumerate}
\usepackage[numbers]{natbib}
\usepackage{subfigure}
\usepackage{setspace}
\usepackage{color}
\usepackage{dsfont}
\usepackage{natbib}

\def \w {\omega}
\def \t {\theta}
\def \a {\alpha}
\def \b {\beta}

\def \g {\gamma}
\def \l {\lambda}
\def \e {\varepsilon}
\def \<{\langle}
\def \>{\rangle}
\def \d {\delta}
\def \O {\Omega}
\def \ad{{\rm ad}}

\def \F{\mathcal{F}}
\def \B{\mathcal{B}}

\def \G{\mathcal{G}}
\def \Lie{{\rm Lie}}

\def \Lie{{\rm Lie}}
\title{\LARGE \bf
Ensemble Control on Lie Groups
}

\author{Wei Zhang and Jr-Shin~Li

\thanks{*This work was supported in part by the National Science Foundation under the award ECCS-1810202 and by the Air Force Office of Scientific Research under the award FA9550-17-1-0166.}
\thanks{W. Zhang is with the Department of Electrical and Systems Engineering, Washington University,
		St. Louis, MO 63130, USA
        {\tt\small wei.zhang@wustl.edu}}%
\thanks{J.-S. Li is with the Department of Electrical and Systems Engineering, Washington University,
        St. Louis, MO 63130, USA
        {\tt\small jsli@wustl.edu}. Questions, comments, or corrections to this document may be directed to J.-S. Li at this email address.}%
}

\begin{document}

\maketitle
\begin{abstract}
Problems involving control of large ensmebles of structurally identical dynamical systems, called \emph{ensemble control}, arise in numerous scientific areas from quantum control and robotics to brain medicine. In many of such applications, control can only be implemented at the population level, i.e., through broadcasting an input signal to all the systems in the population, and this new control paradigm challenges the classical systems theory. In recent years, considerable efforts have been made to investigate controllability properties of ensemble systems, and most works emphasized on linear and some forms of bilinear and nonlinear ensemble systems. In this paper, we study controllability of a broad class of bilinear ensemble systems defined on semisimple Lie groups, for which we define the notion of ensemble controllability through a Riemannian structure of the state space Lie group. Leveraging the Cartan decomposition of semisimple Lie algebras in representation theory, we develop a \emph{covering method} that decomposes the state space Lie group into a collection of Lie subgroups generating the Lie group, which enables the determination of ensemble controllability by controllability of the subsystems evolving on these Lie subgroups. Using the covering method, we show the equivalence between ensemble and classical controllability, i.e., controllability of each individual system in the ensemble implies ensemble controllability, for bilinear ensemble systems evolving on semisimple Lie groups. This equivalence makes the examination of controllability for infinite-dimensional ensemble systems as tractable as for a finite-dimensional single system. 
\end{abstract}

\section{Introduction}

Finely manipulating a large ensemble of structurally identical dynamical systems has emerged as an essential demand in diverse areas from quantum science and technology \cite{Glaser98,Li_PRA06,Li_PNAS11,Dong10,Dong12,Augier18}, brain medicine \cite{Zlotnik12,Ching2013b,Kafashan2015,Li_NatureComm16} and robotics \cite{Becker12} to sociology \cite{Brockett10,Chen16_flocks}. In many applications involving ensemble systems, control can only be exerted at the population level becasue it is infeasible and often impossible to receive state feedback for each individual system. As a result, considerable efforts have been made over the past years to understand the fundamental limit on the extent to which an ensemble system can be manipulated with a broadcast open-loop signal. This new control paradigm raised significant challenges in classical systems theory, while offering abundant opportunities for making theoretical advancements.

Among the developments in this rising area, referred to as ensemble control, extensive focuses have been placed on investigating the controllability property of ensemble systems, including linear \cite{Li_TAC11,Helmke14,Zeng_16_moment,Li_TAC16,Dirr18,Li_SICON20}, bilinear \cite{Li_TAC09,Beauchard10,Chen2020}, and some forms of nonlinear ensemble systems \cite{Li_TAC13,XD_Chen19,Kuritz19}. The work on analyzing  controllability of an ensemble with each system defining on the Lie group SO$(3)$ set the milestone in formal and rigorous study of ensemble systems \cite{Li_TAC09}. In this work, using Lie algebraic tools, the controllability analysis was translated to the problem of polynomial approximation, which opened the door for addressing ensemble control problems from the perspective of ``approximation''. This new notion has led to seminal works on developing necessary and/or sufficient conditions for ensemble controllability \cite{Li_TAC11,Dirr18,Helmke14,Li_TAC16,Zeng_16_moment,Zhang18,Li_SICON20} and observability \cite{zeng2016tac,Zeng16}, and novel theory- and computational-based techniques for optimal ensemble control design and synthesis \cite{Li_PNAS11,Li_ACC12_SVD,Dong14,Gong16,Li_SICON17,Li_Automatica18}. Although progress in understanding fundamental properties of nonlinear ensemble systems is underdeveloped \cite{Li_TAC13,XD_Chen19} and much is awaiting to be explored, the work presented in \cite{Li_TAC09} shed light on revealing the equivalence between ensemble controllability and classical controllability for certain classes of ensemble systems. 

In general, controllability of each individual system (i.e., classical controllability) in an ensemble is a necessary condition to ensemble controllability but not sufficient. Namely, if an ensemble system is ensemble controllable, then each individual system in the ensemble must be controllable in the classical sense; however, the reversal is generally not true. Motivated by the work on the control of ensemble systems on SO(3) \cite{Li_TAC09}, where controllability of each individual system led to controllability of the entire ensemble, in this paper, we extend this previous finding to explore such equivalence in classical and ensemble controllability for more general classes of ensemble systems. Specifically, we study the bilinear ensemble system in which each individual system evolves on the same semisimple Lie group. In our approach, such an ensemble is regarded as a single system defined on the space of Lie group-valued functions, which is an infinite-dimensional Lie group, and the concept of ensemble controllability is rigorously defined in the sense of approximate controllability through a bi-invariant metric on this infinite-dimensional Lie group. 
The main tool developed in this work 
is the \emph{covering method}. The central idea of this method is to decompose the state space Lie group of a bilinear ensemble system into a collection of Lie subgroups, which generates the Lie group, so that controllability of the ensemble is determined by that of the subsystems evolving on these Lie subgroups.
The covering method is further used to reveal a significant consequence of equivalence between ensemble and classical controllability of bilinear systems defined on semisimple Lie groups, i.e., classical controllability of each individual system in the ensemble implies ensemble controllability. Moreover, we show that this equivalence 
is not constrained to systems evolving on compact Lie groups and holds for bilinear ensemble systems induced by Lie group actions on vector spaces, for which each individual system is defined on a non-compact Lie group. 

This paper is organized as follows. In the next section, we introduce the notion of ensemble controllability for parameterized families of control systems evolving on Lie groups through the bi-invariant Riemannian structures of the groups. In Section \ref{sec:so(3)}, we revisit and extend our previous results in ensemble controllability of bilinear systems on SO$(3)$, which lays a foundation for the investigation into controllability of bilinear ensemble systems on general semisimple Lie groups. In Section \ref{sec:ensemble_son}, we introduce the covering method to establish the equivalence between ensemble and classical controllability for bilinear systems. In particular, we first illustrate the main idea by using systems evolving on SO$(n)$ with $n>3$, and then extend the analysis to systems defined on general semisimple Lie groups by using Cartan decompositions. The generality of the equivalence to ensemble systems induced by Lie group actions on vector spaces is presented in Section \ref{sec:sen}.

\section{Preliminaries}
\label{sec:prelim}
In this section, we review the classical controllability results characterized by the Lie algebra rank condition (LARC) for control systems defined on compact, connected Lie groups. Then, we introduce the notion of ensemble controllability for a parameterized family of systems defined on a Lie group through the Riemannian structure of this group, and address the major obstacle to ensemble controllability analysis of such systems when applying LARC.

\subsection{Controllability of systems on compact and connected Lie groups}
Controllability of systems evolving on compact, connected Lie groups has been extensively studied \cite{Brockett72, Jurdjevic72, Jurdjevic96}. 
The analysis is based on examining whether the Lie algebra generated by the drift and control vector fields is equivalent to the underlying Lie algebra of the Lie group. Specifically, a right-invariant bilinear control system defined on a compact, connected Lie group $G$ of the form,
\begin{align}
	\label{eq:general_form}
	\frac{d}{dt}X(t)=\Big[B_0+\sum_{i=1}^mu_i(t)B_i\Big]X(t),\quad X(0)=I,
\end{align}
is of great theoretical and practical interest, where $X(t)\in G$ is the state, $B_0,\dots,B_m$ are elements in the Lie algebra $\mathfrak{g}$ of $G$, 
$I$ is the identity element of $G$, and $u_i(t)\in\mathbb{R}$ are piecewise constant control functions for $i=1,\ldots,m$. In addition, we denote the Lie algebra generated by the set of vector fields $\F=\{B_0,B_1,\dots,B_m\}$ by Lie$\{B_0,B_1,\dots,B_m\}$, i.e., the smallest linear subspace of $\mathfrak{g}$, which contains $\F$ 
and is closed under the Lie bracket operation defined by $[M,N]=MN-NM$ for all $M,N\in\mathfrak{g}$. Controllability of the system of the form in \eqref{eq:general_form} can be evaluated by the following theorem.

\begin{theorem}
	\label{thm:LARC}
	The system in \eqref{eq:general_form} is controllable on the 
Lie group $G$ if and only if $\Lie(\mathcal{F})=\mathfrak{g}$, where $\mathcal{F}=\{B_0,B_1,\dots,B_m\}$. 
\end{theorem}

{\it Proof.}
See \cite{Brockett72, Khaneja00}. \hfill$\Box$

If the dimension of $\mathfrak{g}$ is $n$, then the only linear subspace of $\mathfrak{g}$ that also has dimension $n$ is $\mathfrak{g}$ itself. Thus, checking controllability of a control system as in 
\eqref{eq:general_form} is equivalent to checking the dimension of $\Lie(\F)$. 
Conventionally, the necessary and sufficient condition in Theorem \ref{thm:LARC} is referred to as the Lie algebra rank condition (LARC).

\subsection{Control of ensemble systems}
\label{sec:ensemble_def}
An ensemble control system is a family of control systems defined on a manifold $M$,
\begin{align}
	\label{eq:ensemble}
	\frac{d}{dt}x(t,\b)=f(t,x(t,\b),u(t)),
\end{align}
parameterized by a parameter $\b\in K\subset\mathbb{R}^d$ such that $x(t,\b)\in M$ for each $t\in\mathbb{R}$ and $\b\in K$, 
where the parameter space $K$ is generally assumed to be compact. In this case, for each fixed $t\in\mathbb{R}$, $x(t,\cdot)$ is an $M$-valued function defined on $K$, i.e., the state space of the ensemble system in \eqref{eq:ensemble} is actually a space of $M$-valued functions defined on $K$, 
denoted by $\F(K,M)$. The parameter independent open-loop control input $u(t)\in\mathbb{R}^m$ is a broadcast signal that simultaneously manipulates the ensemble between desired functions in $\F(K,M)$. Note that when the parameter space $K$ is an infinite set, i.e., the ensemble system in \eqref{eq:ensemble} contains infinitely many dynamic units, $\F(K,M)$ is an infinite-dimensional manifold so that the ensemble system is an infinite-dimensional system. For such systems, we define the notion of ensemble controllability in the approximation sense.

\begin{definition}[Ensemble Controllability]
	\label{def:ec}
	Let $\mathcal{F}(K,M)$ denote a space of $M$-valued functions defined on $K$. The family of systems in \eqref{eq:ensemble} is said to be ensemble controllable on the function space $\mathcal{F}(K,M)$, if for any $\e>0$ and starting with any initial state $x_0\in\mathcal{F}(K,M)$, where $x_0(\cdot)=x(0,\cdot)$, there exists a control law $u(t)$ that steers the system into an $\e$-neighborhood of a desired target state $x_F\in\mathcal{F}(K,M)$ at a finite time $T>0$, i.e., $d(x(T,\cdot),x_F(\cdot))<\e$, where $d:\F(K,M)\times\F(K,M)\rightarrow\mathbb{R}$ is a metric on $\mathcal{F}(K,M)$. Note that the final time $T$ may depend on $\e$, and ensemble controllability is a notion of approximate controllability.
\end{definition}

In this work, we focus on the time-invariant bilinear ensemble system evolving on a Lie group $G$ of the form
\begin{align}
	\label{eq:ensemble_Lie}
	\frac{d}{dt}X(t,\b)=\Big[\b_0B_0+\sum_{i=1}^m\b_i \, u_i(t)B_i\Big]X(t,\b), \quad X(0,\b)=I,
\end{align}
where $\b=(\b_0,\dots,\b_m)'$ is the parameter vector varying on a compact subset $K\subset\mathbb{R}^{m+1}$, $X(t,\cdot)\in C(K,G)$ is the state and $C(K,G)$ denotes the space of continuous $G$-valued functions defined on $K$, $B_0,\dots,B_m$ are elements in the Lie algebra $\mathfrak{g}$ of $G$,
$I$ is the identity element of $G$, and $u_1,\dots,u_m$ are real-valued piecewise constant control inputs. 

According to Definition \ref{def:ec}, a metric on $C(K,G)$ is necessary in the study of ensemble controllability of the system in \eqref{eq:ensemble_Lie}. In the next section, we will introduce metrics on $C(K,G)$ and $C(K,\mathfrak{g})$ through a Riemannian structure of $G$ such that these two metrics are locally compatible with respect to the exponential map, $\exp:\mathfrak{g}\rightarrow G$. Consequently, ensemble controllability of systems defined on $C(K,G)$ can be studied through their drift and control vector fields in $C(K,\mathfrak{g})$.


\subsection{Metric space structures on $C(K,G)$}
\label{sec:metric_G}
In Definition \ref{def:ec}, ensemble controllability is defined in the sense of approximate controllability, where it only requires to steer the considered system into an $\e$-neighborhood of the desired final state. However, the properties of neighborhoods depend on the topology of the state space of the system. Therefore, in this section, we will introduce a metrizable topology on $C(K,G)$ such that ensemble controllability of an ensemble system evolving on $C(K,G)$ can be defined through the metric induced by this topology.

The compact-open topology is commonly used on the space of continuous functions between two topological spaces. In our case, $K$ is compact and $G$ is a metric space as a Riemannian manifold, then the compact-open topology on $C(K,G)$ is metrizable. Specifically, it is equivalent to the topology of uniform convergence \cite{Hatcher02}, i.e., the topology induced by the metric $d(f,g)=\sup_{\b\in K}\rho(f(\b),g(\b))$ for any $f,g\in C(K,G)$, where $\rho:G\times G\rightarrow G$ is the metric induced by a Riemannian metric on $G$. This observation illustrates that it suffices to define a Reimannian structure on $G$, which in turn induces a metric on $C(K,G)$.

A bi-invariant Riemannian metric is a good candidate of Reimannian metrics defined 
on a compact, connected Lie group $G$ for understanding the relationship between its geometric and algebraic structures. Because, under this metric, the exponential map from $\mathfrak{g}$ to $G$ coincides with the Riemannian exponential map from $T_IG$ to $G$, where $T_IG$ denotes the tangent space of $G$ at the identity element $I$ \cite{Petersen16}. Correspondingly, the trajectory of each individual system in the ensemble in \eqref{eq:ensemble_Lie} 
is a concatenation of some geodesics of $G$. Computationally, a bi-invariant Riemannian metric can be obtained by averaging an arbitrary inner product defined on $\mathfrak{g}$ over the group $G$, where $\mathfrak{g}$ is identified with $T_IG$ of $G$ \cite{Warner83}. 


Let $\<\cdot,\cdot\>:\mathfrak{g}\times\mathfrak{g}\rightarrow\mathbb{R}$ denote an inner product on $\mathfrak{g}$ that extends to a bi-invariant metric on $G$, then the uniform norm on $C(K,\mathfrak{so}(n))$, that is, $\|f-g\|_{\infty}=\sup_{\b\in K}\|f(\b)-g(\b)\|$ for $f,g\in C(K,\mathfrak{so}(n))$, is well-defined because $K$ is compact, where $\|f(\beta)-g(\beta)\|=\sqrt{\<f(\b)-g(\b),f(\b)-g(\b)\>}$ is the norm on $\mathfrak{g}$ induced by the inner product. If $\|f-g\|_{\infty}<\e$ for some $\e$ smaller than the injectivity radius of the Riemannian exponential map, then $\rho(\exp(f(\beta)),\exp(g(\beta)))\leq\|f(\b)-g(\b)\|\leq\|f-g\|_{\infty}<\e$ holds for any $\b\in K$, because the Lie group $G$ with the bi-invariant Riemannian metric has non-neagtive sectional curvature \cite{Petersen16}, where $\rho$ is the metric on $G$ induced by the bi-invariant Riemannian metric. On the other hand, since 
$G$ is connected and compact, the exponential map $\exp:\mathfrak{g}\rightarrow G$ is surjective \cite{Hall15}, and thus the uniform topology of $C(K,G)$ is carried over from the uniform norm of $C(K,\mathfrak{g})$. This property 
enables the study of ensemble controllability of the system in \eqref{eq:ensemble_Lie} on $C(K,G)$ through its drift and control vector fields on $C(K,\mathfrak{g})$.

It can be shown that $C(K,G)$ itself is an infinite-dimensional Lie group with the Lie algebra $C(K,\mathfrak{g})$. Furthermore, since every element $f\in C(K,\mathfrak{g})$ can be expressed in the form $f=\sum_{i=1}^n f_i E_i$ for some $f_i\in C(K,\mathbb{R})$ with $\{E_1,\dots, E_n\}$ a basis of $\mathfrak{g}$, this indicates that $C(K,\mathfrak{g})$, as a $C(K,\mathbb{R})$-module, is isomorphic to $C(K,\mathbb{R})\otimes\mathfrak{g}$,
where $C(K,\mathbb{R})$ is the set of continuous real-valued functions defined on $K$ and $\otimes$ denotes the tensor product over $\mathbb{R}$. 
However, $C(K,\mathbb{R})$ is generally not compact with respect to the topology of uniform convergence, e.g., the sequence $f_n(\b)=\b^n$ in $C([0,1],\mathbb{R})$ has no convergent subsequence. Consequently, $C(K,G)$ is a non-compact infinite-dimensional Lie group, which disables 
the application of the LARC, as presented in Theorem \ref{thm:LARC}, to examine controllability of ensemble systems defined on $C(K,G)$ and hence motivates the need of developing new tools to achieve this goal. 

To this end, in Sections \ref{sec:so(3)} and \ref{sec:ensemble_son}, we integrate tools from geometry, analysis, and algebra to synthesize the machinery for controllability analysis of ensemble systems defined on $C(K,G)$ in the form of \eqref{eq:ensemble_Lie}. In particular, our framework will be elaborated through 
the study of the ensemble system defined on $C(K,{\rm SO}(n))$ by leveraging the nice structure of $\mathfrak{so}(n)$, where SO$(n)$ is the special orthogonal group consisting of all $n$-by-$n$ orthogonal matrices with determinant 1 and $\mathfrak{so}(n)$ is its Lie algebra consisting of all $n$-by-$n$ skew-symmetric matrices. In the next section, we will initiate our investigation with the ensemble system evolving on $C(K,{\rm SO}(3))$.

\section{Ensemble control of systems on SO(3)} 
\label{sec:so(3)}
Manipulating an ensemble of systems evolving on SO(3) is an important problem arising in many areas, notably in quantum control and robotics \cite{Glaser98,Li_PRA06,Dong10,Dong12,Li_NatureComm17,Becker12}. 
In this section, we revisit and extend our previous results in ensemble controllability of systems on SO(3) \cite{Li_TAC09}, which will lay the foundation for analyzing controllability of ensemble systems defined on SO$(n)$ and, further, on SE$(n)$.

We first consider the driftless ensemble system on SO$(3)$, given by
\begin{align}
	\label{eq:ensemble_SO3}
	\frac{d}{dt}X(t,\b)=\b\big[u\O_y+v\O_x\big]X(t,\b), \quad X(0,\b)=I, 
\end{align}
where
$\b\in K=[a,b]\subset\mathbb{H}$, $\mathbb{H}=\mathbb{R}^+=(0,\infty)$, and
\begin{align*}
	\O_y=\left[\begin{array}{ccc} 0 & 0 & 1 \\ 0 & 0 & 0 \\ -1 & 0 & 0 \end{array}\right], \quad \O_x=\left[\begin{array}{ccc} 0 & 0 & 0 \\ 0 & 0 & -1 \\ 0 & 1 & 0 \end{array}\right]
\end{align*}
are the generators of rotation around the $y$- and the $x$-axis, respectively. According to the discussion in Section \ref{sec:prelim}, a metric on $C(K,{\rm SO}(3))$ is required to define the notion of ensemble controllability for the system in \eqref{eq:ensemble_SO3}. The detailed construction of a bi-invariant metric on $C(K,{\rm SO}(n))$ is shown in Section \ref{sec:metric_SO(n)}. At present, let's assume that the state space $C(K,{\rm SO}(3))$ has already been equipped with a bi-invariant metric $d:C(K,{\rm SO}(3))\times C(K,{\rm SO}(3))\rightarrow \mathbb{R}$, which is induced by an inner product on $\mathfrak{so}(3)$.
Then, in the following lemma, we prove ensemble controllability of the system in \eqref{eq:ensemble_SO3} over the topology induced by $d$.

\begin{lemma}
	\label{lem:so3}
	The system in \eqref{eq:ensemble_SO3} is ensemble controllable on $C(K,{\rm SO}(3))$.
\end{lemma}

{\it Proof.} We revisit the proof in our previous work \cite{Li_TAC09} by using the metric space structure on $C(K,{\rm SO}(3))$ introduced above. Observe that the Lie brackets generated by the set of matrices $\{\b\O_y,\b\O_x\}$ are
\begin{align*}
	\ad_{\b\O_y}^{2k+1}(\b\O_x) &= (-1)^k\b^{2k}\O_z,\\
	 \ad_{\b\O_y}^{2k}(\b\O_x) &= (-1)^k\b^{2k+1}\O_x,
\end{align*}
where $\ad_AB=[A,B]$ and $\ad_A^kB=[A,\ad_A^{k-1}B]$, $k\in\mathbb{N}$, for all $A,B\in\mathfrak{so}(3)$, and
$$\O_z=\left[\begin{array}{ccc} 0 & -1 & 0 \\ 1 & 0 & 0 \\ 0 & 0 & 0 \end{array}\right]$$
is the generator of rotation around the $z$-axis. Now using elements in $\{\b\O_x,$ $\b^3\O_x,$ $\dots,$ $\b^{2n+1}\O_x\}$ as generators, we are able to produce an evolution of the form
\begin{align}
	R_x(\b) &= \exp(c_0\b\O_x)\exp(c_1\b^3\O_x)\cdots\exp(c_n\b^{2n+1}\O_x) \nonumber \\
	\label{eq:rotation}
	&= \exp\Big\{\sum_{k=0}^nc_k\b^{2k+1}\O_x\Big\}\doteq\exp\left\{\hat{\theta}_x(\b)\O_x\right\}.
\end{align}
As a result, given any $\b$-dependent rotation $\exp\{\t_x(\b)\O_x\}$ around $x$-axis with $\t_x\in C(K,\mathbb{R})$, the order of the polynomial $n$ and the coefficients $c_k$ can be appropriately chosen so that 
$\|\hat\theta_x-\theta_x\|_{\infty}=\sup_{\beta\in K}\sqrt{\<\hat\theta_x(\beta)-\theta_x(\beta),\hat\theta_x(\beta)-\theta_x(\beta)\>}<\e$ for any given approximation error $\e>0$ by the Weierstrass theorem \cite{Apostol74}. Similar arguments can be developed to show that any $\b$-dependent rotations $\exp\{\theta_y(\b)\O_y\}$ and $\exp\{\theta_z(\b)\O_z\}$ around the $y$- and the $z$-axis, respectively, can be approximately generated as $\exp\{\hat{\theta}_y(\b)\O_y\}$ and $\exp\{\hat{\theta}_z(\b)\O_z\}$, and hence any three-dimensional rotations can also be uniformly approximated.  Namely, given any $\beta$-dependent rotation $\Theta\in C(K,{\rm SO}(3))$, 
one can parameterize it by using the Euler angles $\Theta=(\t_x,\t_y,\t_z)$ such that
\begin{align*}
\Theta(\b)&=\exp\{\t_x(\b)\O_x\}\exp\{\t_y(\b)\O_y\}\exp\{\t_z(\b)\O_z\}\\
&=\Theta_x(\b)\Theta_y(\b)\Theta_z(\b),
\end{align*}
and then the desired rotation $\Theta(\b)$ characterized by the three continuous functions, $\t_x,\t_y,\t_z\in C(K,\mathbb{R})$, can be synthesized by using piecewise constant control vector fields as described in \eqref{eq:rotation}. Specifically, for any $\e>0$, the approximated rotations $\hat\t_x$, $\hat\t_y$, and $\hat\t_z$ can be generated such that $\|\hat\t_x-\t_x\|_{\infty}<\e/3$, $\|\hat\t_y-\t_z\|_{\infty}<\e/3$, and $\|\hat\t_z-\t_z\|_{\infty}<\e/3$. As a result, the total evolution
\begin{align*}
\widehat{\Theta}(\b)&=\exp\{\hat\t_x(\b)\O_x\}\exp\{\hat\t_y(\b)\O_y\}\exp\{\hat\t_z(\b)\O_z\}\\
&=\widehat{\Theta}_x(\b)\widehat{\Theta}_y(\b)\widehat{\Theta}_z(\b)
\end{align*}
satisfies
\begin{align*}
d(\widehat{\Theta},\Theta)&=d(\widehat{\Theta}_x\widehat{\Theta}_y\widehat{\Theta}_z,\Theta_x\Theta_y\Theta_z)\\
&\leq d(\widehat{\Theta}_x\widehat{\Theta}_y\widehat{\Theta}_z,\Theta_x\Theta_y\widehat{\Theta}_z)+d(\Theta_x\Theta_y\widehat{\Theta}_z,\Theta_x\Theta_y\Theta_z)\\
&=d(\widehat{\Theta}_x\widehat{\Theta}_y,\Theta_x\Theta_y)+d(\widehat{\Theta}_z,\Theta_z)\\
&\leq d(\widehat{\Theta}_x\widehat{\Theta}_y,\Theta_x\widehat{\Theta}_y)+d(\Theta_x\widehat{\Theta}_y,\Theta_x\Theta_y)+d(\widehat{\Theta}_z,\Theta_z)\\
&=d(\widehat{\Theta}_x,\Theta_x)+d(\widehat{\Theta}_y,\Theta_y)+d(\widehat{\Theta}_z,\Theta_z)\\
&\leq\|\hat\t_x-\t_x\|_{\infty}+\|\hat\t_y-\t_y\|_{\infty}+\|\hat\t_z-\t_z\|_{\infty}<\e,
\end{align*}
where we repeatedly used the triangle inequality and bi-invariance of the metric $d$. This then concludes ensemble controllability of the system in \eqref{eq:ensemble_SO3} on $C(K, {\rm SO}(3))$. \hfill$\Box$

\begin{remark}[Topological characterization of ensemble controllability]
\label{rmk:ensemble_control}
\rm
In the proof of Lemma \ref{lem:so3}, the key observation leading to ensemble controllability of the system in \eqref{eq:ensemble_SO3} is the uniform approximation of $\b$-dependent rotations $\t_x(\b)\O_x,$ $\t_y(\b)\O_y,$ and $\t_z(\b)\O_z$ by iterated Lie bracketing the control vector fields in $\mathcal{G}=\{\b\O_{x},\b\O_y\}$. This implies that the closure of the Lie algebra generated by $\mathcal{G}$ satisfies $\overline{\Lie(\mathcal{G})}=C(K,\mathbb{R})\otimes \mathfrak{so}(3)=C(K,\mathfrak{so}(3))$, which gives rise to a topological characterization of ensemble controllability of the system in \eqref{eq:ensemble_SO3} on $C(K,{\rm SO}(3))$. 
In general, a family of driftless bilinear systems defined on a compact, connected Lie group $G$ parameterized by a vector $\b=(\b_0,\dots,\b_m)'$ varying on a compact subset $K\subset\mathbb{R}^m$ of the form
$$\frac{d}{dt}X(t,\b)=\Big[\sum_{i=1}^m\b_i \, u_i(t)B_i\Big]X(t,\b), \quad X(0,\b)=I,$$
is ensemble controllable on $C(K,G)$ if and only if $\overline{{\rm Lie}(\mathcal{G})}=C(K,\mathfrak{g})$, where $\mathcal{G}$ $=$ $\{\b_1B_1,$ $\dots,$ $\b_m B_m\}$ is the set of control vector fields evaluated at the identity element $I$ of $G$, and $\mathfrak{g}$ is the Lie algebra of $G$. 
\end{remark}


It was also shown in our previous work that the ensemble with a dispersion in the drift, i.e., the system
\begin{equation*}
	\frac{d}{dt}X(t,\b,\w)=\big[\w\O_z+ \b u\O_y+\b v\O_z\big]X(t,\b,\w),\quad X(0,\beta,\omega)=I,
\end{equation*}
where $\w\in K_d\subset\mathbb{R}$ with $K_d$ compact, is ensemble controllable on $C(K\times K_d,{\rm SO}(3))$ \cite{Li_TAC09}. In the following, we illustrate the applicability of the polynomial approximation technique exploited in the proof of Lemma \ref{lem:so3} to analyze ensemble systems on SO(3) with three parameter variations. This analysis constitutes the key element in the covering method to be developed in Section \ref{sec:ensemble_son} for the controllability analysis of bilinear ensemble systems defined on compact, connected Lie groups.


\begin{proposition}
	\label{prop:so3_3dispersion}
	An ensemble system of the form,
	\begin{align}
		\label{eq:so3_3dispersion}
		\frac{d}{dt}X(t,\b)=\big[\b_1u_1\O_x+\b_2u_2\O_y+\b_3u_3\O_z\big]X(t,\b),\quad X(0,\b)=I,
	\end{align}
	is ensemble controllable on $C(K,{\rm SO}(3))$, where $\b=(\b_1,\b_2,\b_3)\in K$ is the parameter vector varying on a compact subset $K$ of the three-dimensional upper half space $\mathbb{H}^3=\{(\beta_1,\beta_3,\beta_3)\in\mathbb{R}^3:\beta_i>0\text{ for all }i=1,\dots,3\}$, $I$ is the 3-by-3 identity matrix, and $u_i(t)$ are piecewise constant control inputs for all $i=1,2,3$.	
\end{proposition}

{\it Proof.}
	By successive Lie brackets of the control vector fields $\b_2\O_y$ and $\b_3\O_z$, we obtain
\begin{align*}
	\ad_{\b_2\O_y}^{2k+1}(\b_3\O_z) &= (-1)^k\b_2^{2k+1}\b_3\O_x, \\
	\ad_{\b_3\O_z}^{2l+1}(\b_2^{2k+1}\b_3\O_x) &= (-1)^{l}\b_2^{2k+1}\b_3^{2l+1}\O_x,
\end{align*}
where $k,l\in\mathbb{N}$. Then, defining $L_{(k,l)}=\b_2^{2k+1}\b_3^{2l+1}$ and applying iterated Lie brackets of $[\b_1\O_x,\b_2\O_y]$ and $L_{(k,l)}\O_x$ yields
\begin{align*}
	\ad_{[\b_1\O_x,\b_2\O_y]}^{2s}(L_{(k,l)}\O_x) &= (-1)^s\b_1^{2s}\b_2^{2(k+s)+1}\b_3^{2l+1}\O_x\\
	&= (-1)^s\b_1^{2s}\b_2^{2(k+s)}\b_3^{2l}(\b_2\b_3\O_x),
\end{align*}
where $s\in\mathbb{N}$. Furthermore, let $L_{(s,k,l)}(\b)=\b_1^{2s}\b_2^{2(k+s)}\b_3^{2l}$ and $\mathcal{A}={\rm span}\{L_{(s,k,l)}:s,k,l=0,1,\dots\}\subset C(K,\mathbb{R})$, then we claim that $\mathcal{A}$ is a subalgebra of $C(K,\mathbb{R})$ by checking that $fg\in\mathcal{A}$ for any $f,g\in\mathcal{A}$. Now, pick any two points $x=(x_1,x_2,x_3)'$ and $y=(y_1,y_2,y_3)'$ in $K$ and assume $f(x)=f(y)$ for all $f\in\mathcal{A}$, in particular,  $L_{(1,0,0)}(x)=L_{(1,0,0)}(y)$, $L_{(0,1,0)}(x)=L_{(0,1,0)}(y)$, and $L_{(0,0,1)}(x)=L_{(0,0,1)}(y)$ hold. This gives $x_i=y_i$ for each $i=1,2,3$, i.e., $x=y$. Therefore, $\mathcal{A}$ separates points in $K$ \cite{Li_SICON20} 
and hence $\mathcal{A}$ is dense in $C(K,\mathbb{R})$ by Stone-Weierstrass Theorem \cite{Folland99}. 
Equivalently, for any $f\in C(K,\mathbb{R})$, we can uniformly approximate $f(\b)\O_x$ by iterated Lie brackets of the control vector fields in $\mathcal{G}=\{\b_1\O_x,\b_2\O_y,\b_3\O_z\}$. A similar argument can be applied to show that, for any $g,h\in C(K,\mathbb{R})$, $g(\b)\O_y$ and $h(\b)\O_z$ can also be uniformly approximated. It follows that $\overline{{\rm Lie}(\mathcal{G})}=C(K,\mathbb{R})\otimes\mathfrak{so}(3)=C(K,\mathfrak{so}(3))$, 
and hence the system in \eqref{eq:so3_3dispersion} is ensemble controllable on $C(K,\rm{SO(3)})$ by Remark \ref{rmk:ensemble_control}.\hfill$\Box$

\section{Ensemble control of systems on compact Lie groups}
\label{sec:ensemble_son}
In this section, we will carry out an extension of the ensemble controllability analysis developed in the previous section dedicated to the system on SO$(3)$ to general systems defined on compact, connected Lie groups. To this end, we will introduce a \emph{covering method} based on the decomposition of the state space Lie group into a collection of Lie subgroups, which generates this Lie group,  and, correspondingly, decomposes the ensemble system defined on this Lie group into a collection of subsystems, each of which evolves on one of these Lie subgroups. This decomposition then enables the determination of controllability of the ensemble by controllability of each subsystem, since the state space Lie group is generated by the Lie subgroups defining the state space of the subsystems.


Before the discussion of systems evolving on general semisimple Lie groups, this method will be best motivated and illuminated with the system defined on SO$(n)$ first. To facilitate our exposition, we review some key properties of the Lie algebra $\mathfrak{so}(n)$ that are relevant to the subsequent ensemble controllability analysis in the following section. 

\subsection{Basics of the Lie algebra $\mathfrak{so}(n)$}
\label{sec:basics}
The Lie algebra $\mathfrak{so}(n)$ is the vector space containing all $n\times n$ real skew-symmetric matrices, which has dimension $n(n-1)/2$. Let $E_{ij}\in\mathbb{R}^{n\times n}$ denote the matrix whose $ij^{\rm th}$ entry is $1$ and others are $0$, then the matrix $\O_{ij}=E_{ij}-E_{ji}$ satisfies
\begin{align*}
	\O_{ij}=
	\begin{cases}
	-\O_{ji},\ \text{if\ } i\neq j,\\
	0, \qquad \text{if\ } i=j,
	\end{cases}
\end{align*}
taking value 1 in the $ij^{\rm th}$ entry, -1 in the $ji^{\rm th}$ entry, and 0 elsewhere. Moreover, the set $\mathcal{B}=\{\O_{ij}:1\leq i<j\leq n\}$ forms a basis of $\mathfrak{so}(n)$, which is referred to as the \emph{standard basis} of $\mathfrak{so}(n)$. 

\begin{lemma}
	\label{lem:son}
	The Lie bracket of $\O_{ij}$ and $\O_{kl}$ 
	satisfies the relation $[\O_{ij},\O_{kl}]=\d_{jk}\O_{il}+\d_{il}\O_{jk}+\d_{jl}\O_{ki}+\d_{ik}\O_{lj}$, where $\d$ is the Kronecker delta function, i.e.,
	\begin{align*}
		\d_{mn}=\begin{cases} 1, \quad {\rm if\ } m=n, \\ 0, \quad {\rm if\ } m\neq n. \end{cases}
	\end{align*}
\end{lemma}
{\it Proof.}
Notice that $E_{ij}E_{kl}=\d_{jk}E_{il}$, so $[E_{ij},E_{kl}]=\delta_{jk}E_{il}-\delta_{li}E_{kj}$. Following the bilinearity of the Lie bracket, we get
\begin{align*}
	[\O_{ij},\O_{kl}]&=[E_{ij}-E_{ji},E_{kl}-E_{lk}]= [E_{ij},E_{kl}]-[E_{ij},E_{lk}]-[E_{ji},E_{kl}]+[E_{ji},E_{lk}] \\
	&= \delta_{jk}E_{il}-\delta_{li}E_{kj}-\delta_{jl}E_{ik}+\delta_{ki}E_{lj}-\delta_{ik}E_{jl}+\delta_{lj}E_{ki}+\delta_{il}E_{jk}-\delta_{kj}E_{li}\\
	&= \delta_{jk}\O_{il}+\delta_{il}\O_{jk}+\delta_{jl}\O_{ki}+\delta_{ik}\O_{lj}. 
\end{align*}\hfill$\Box$

According to Lemma \ref{lem:son}, for any $\O_{ij},\O_{kl}\in\mathcal{B}$, $[\O_{ij},\O_{kl}]\neq 0$ if and only if $i=l$ $j=k$, $i=k$ or $j=l$. 


\subsection{Bi-invariant metrics on ${\rm SO}(n)$}
\label{sec:metric_SO(n)}
By Definition \ref{def:ec} in Section \ref{sec:ensemble_def}, a metric on $C(K,{\rm SO}(n))$ is required to define the notion of ensemble controllability for systems evolving on SO$(n)$. Moreover, because SO$(n)$ is a Lie group, the discussion in Section \ref{sec:metric_G} implies that a metric on $C(K,{\rm SO}(n))$ can be induced by an inner product on the Lie algebra $\mathfrak{so}(n)$. In particular, we introduce an inner product $\<\cdot,\cdot\>:\mathfrak{so}(n)\times\mathfrak{so}(n)\rightarrow\mathbb{R}$ such that the standard basis elements in $\B$ form an orthonormal basis for $\mathfrak{so}(n)$, or equivalently, $\<\O_{ij},\O_{kl}\>={\rm tr}(\O_{ij}'\O_{kl})/2$. Then, we extend this inner product to a left-invariant Riemannian metric on SO$(n)$ by defining $\<\O_{ij}X,\O_{kl}X\>={\rm tr}(\O_{ij}'\O_{kl})/2$ for any $X\in{\rm SO}(n)$. Notice that $\<\cdot,\cdot\>$ is invariant under the adjoint action of SO$(n)$ on $\mathfrak{so}(n)$, i.e., $\<XYX^{-1},XZX^{-1}\>=\<Y,Z\>$ for any $X\in{\rm SO}(n)$ and $Y,Z\in\mathfrak{so}(n)$. Hence, this left-invariant Riemannian metric is also bi-invariant \cite{Petersen16}, which then induces a bi-invariant metric $\rho$ on SO$(n)$. Consequently, by the discussion in Section \ref{sec:metric_G}, the compact-open topology induces a bi-invariant metric $d$ on $C(K,{\rm SO}(n))$, which coincides with the topology of uniform convergence with respect to $\rho$, i.e., $d(f,g)=\sup_{\beta\in K}\rho(f(\beta),g(\beta))$ for any $f,g\in C(K,{\rm SO}(n))$. In particular, for the 
case of SO$(3)$ discussed in Section \ref{sec:so(3)}, the bi-invariant metric $d$ is just obtained by defining the set $\{\O_x,\O_y,\O_z\}$ to be an orthonormal basis of $\mathfrak{so}(3)$. 

In the following sections, ensemble controllability will be analyzed under this bi-invariant metric $d$ on $C(K,{\rm SO}(n))$.


\subsection{The covering method for ensemble controllability analysis}
\label{sec:covering}
In this section, we develop a covering method for examining ensemble controllability of bilinear systems evolving on semisimple Lie groups. Together with the technique of polynomial approximation, we then establish an equivalence between ensemble and classical controllability for such bilinear ensemble systems. The existence and construction of this covering method are based on the Cartan decomposition of semisimple Lie algebras in representation theory \cite{Hall15}. Specifically, given such a system, we apply the Cartan decomposition to the semisimple Lie algebra of the state-space Lie group, which gives rise to a cover of the Lie algebra consisting of Lie subalgebras isomorphic to $\mathfrak{so}(3)$ or $\mathfrak{su}(2)$. Correspondingly, the ensemble system also admits a decomposition into a family of ensemble subsystems with each 
defined on SO$(3)$ or SU$(2)$. In this way, the controllability analysis of the ensemble system is equivalently carried over to these ensemble subsystems.
To showcase the main idea of the decomposition in the covering method, we use an example of the Lie group SO$(4)$. 


\begin{example}[A simple illustration of the covering method]\rm
	\label{ex:cover}
	In this example, we will construct a set of generators of SO$(4)$ such that every generator is a Lie subgroup of SO$(4)$ isomorphic to SO$(3)$. We start our construction with decomposing the Lie algebra $\mathfrak{so}(4)$ into a collection of Lie subalgebras isomorphic to $\mathfrak{so}(3)$. This is equivalent to constructing a cover of the standard basis $\B=\{\O_{12},\O_{13},\O_{14},\O_{23},\O_{24},\O_{34}\}$. To this end, let $\mathcal{U}=\{\B_1,\B_2,\B_3,\B_4\}$, where $\B_1=\{\O_{12},\O_{13},\O_{23}\}$, $\B_2=\{\O_{12},\O_{24},\O_{14}\}$, $\B_3=\{\O_{13},\O_{14},\O_{34}\}$, and $\B_4=\{\O_{23},\O_{34},\O_{24}\}$, then it is clear that $\mathcal{U}$ forms a cover of $\mathcal{B}$, because $\B=\B_1\cup\B_2\cup\B_3\cup\B_4$. Moreover, let $F=\{\Lie(\B_1),\Lie(\B_2)$, $\Lie(\B_3),\Lie(\B_4)\}$, then we have ${\rm span}(F)=\mathfrak{so}(4)$, 
	and hence $F$ is a set of generators of $\mathfrak{so}(4)$. Notice that each $\Lie(\B_i)$, $i=1,\ldots,4$, is isomorphic to $\mathfrak{so}(3)$ so that its Lie group  $G_i$ is a Lie subgroup of ${\rm SO}(4)$ isomorphic to SO$(3)$.  In addition, because $F$ 
	generates $\mathfrak{so}(4)$, $\mathcal{V}=\{G_1,G_2,G_3,G_4\}$ is a set of generators of SO$(4)$ as desired. This cover of SO$(4)$ is illustrated in Figure \ref{fig:cover}.

\begin{figure}[t]
     \centering
     \includegraphics[width=0.7\columnwidth]{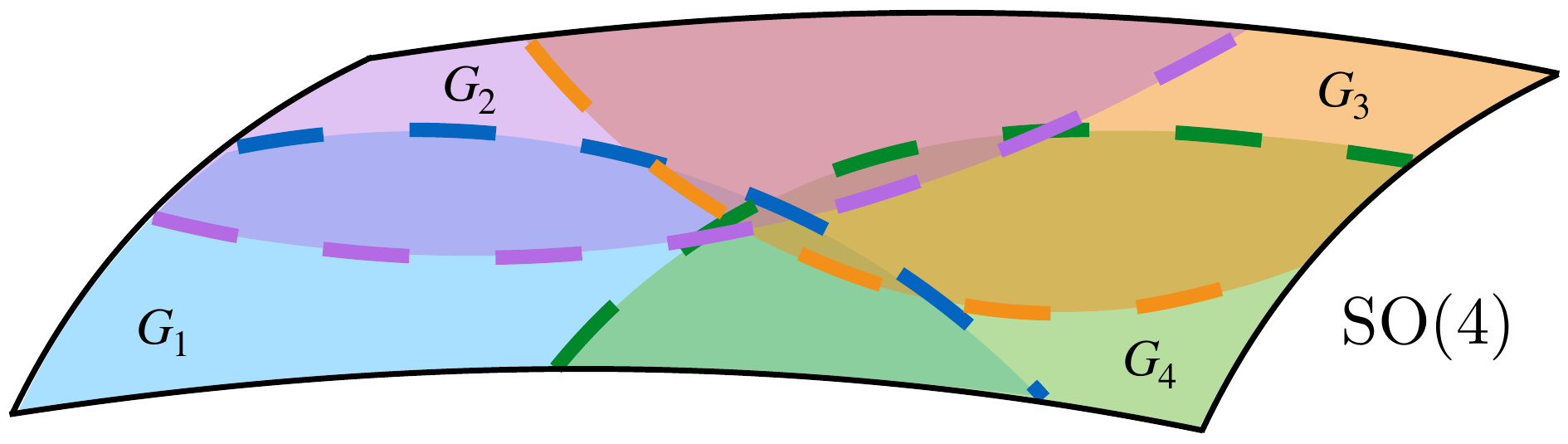}
     \caption{The demonstration of the cover $\mathcal{V}=\{G_1,G_2,G_3,G_4\}$ of {\rm SO}$(4)$ constructed in Example \ref{ex:cover}. In particular, $G_1$, $G_2$, $G_3$, and $G_4$, illustrated by blue, purple, orange, and green shadows bounding by the dashed lines with the corresponding colors, respectively, are Lie subgroups of {\rm SO}$(4)$ isomorphic to {\rm SO}$(3)$.}
     \label{fig:cover}
\end{figure}
\end{example}

The covering idea 
illustrated in Example \ref{ex:cover} for SO$(4)$ can be directly generalized to SO$(n)$. This generalization immediately enables the adoption of the polynomial approximation based technique developed for systems on SO$(3)$ in Section \ref{sec:so(3)} to the ensemble controllability analysis of systems on SO$(n)$ with $n>3$. More importantly, the covering method paves the way for understanding and quantifying the equivalence between ensemble and classical controllability. 


\begin{theorem}[The main result]
	\label{thm:son_ensemble}
	Consider an ensemble of systems on {\rm SO}$(n)$, given by
	\begin{align}
		\label{eq:son_ensemble}
		\frac{d}{dt}X(t,\b)=\Big[\sum_{k=1}^m \b_k u_k(t)\, \Omega_{i_kj_k}\Big]X(t,\b),\quad X(0,\b)=I,
	\end{align}
	where the parameter vector $\b=(\b_1,\dots,\b_m)'$ takes values on a compact subset $K\subset\mathbb{H}^m$, the state $X(t,\cdot)\in C(K,{\rm SO}(n))$, and the control inputs $u_k(t)\in\mathbb{R}$ are piecewise constant for all $k=1,\dots,m$. This system is  ensemble controllable on $C(K,{\rm SO}(n))$ if and only if each individual system with respect to 
	a fixed $\beta\in K$ in this ensemble is controllable on ${\rm SO}(n)$.
\end{theorem}
{\it Proof.} 
The necessity is obvious, and hence it remains to show the sufficiency. In particular, we divide the proof of sufficiency into three steps. 

(Step I): An ensemble of systems defined on SO$(n)$ of the form,
\begin{align}
\label{eq:son_ensemble_full}
	\frac{d}{dt}X(t,\b)=\Big[\sum_{1\leq i<j\leq n}\b_{ij}u_{ij}(t)\O_{ij}\Big]X(t,\b),\quad  X(0,\b)=I,
\end{align}
is ensemble controllable on $C\big(\prod_{1\leq i<j\leq n}K_{ij},{\rm SO}(n)\big)$, where the parameter vector $\b=(\b_{12},\dots,\b_{n-1,n})$ takes values in the product space $\prod_{1\leq i<j\leq n}K_{ij}$ with each $K_{ij}$ a compact subset of $\mathbb{H}$, $X(t,\cdot)\in C\big(\prod_{1\leq i<j\leq n}K_{ij},{\rm SO}(n)\big)$ is the state, and $u_{ij}(t)\in\mathbb{R}$ are piecewise constant for all $1\leq i<j\leq n$. 

For any $\O_{ij}\in\mathcal{B}$ and $k_1\in\{1,\dots,n\}\backslash\{i,j\}$, the subset $\mathcal{S}_1=\{\O_{ij},\O_{ik_1},\O_{k_1j}\}$ of $\mathcal{B}$ generates a Lie subalgebra of $\mathfrak{so}(n)$ isomorphic to $\mathfrak{so}(3)$. By Proposition \ref{prop:so3_3dispersion}, the controllable submanifold of the system obtained by setting $u_{\a\g}=0$ for all $\a,\g\in\{1,\dots,n\}\backslash\{i,j,k_1\}$ in the system \eqref{eq:son_ensemble_full}, i.e., 
\begin{align*}
&\frac{d}{dt}X(t,\b)=[\b_{ij}u_{ij}(t)\O_{ij}+\b_{ik_1}u_{ik_1}(t)\O_{ik_1}+\b_{k_1j}u_{k_1j}(t)\O_{k_1j}]X(t,\beta),\\
&X(0,\b)=I,
\end{align*}
is a Lie subgroup of $C(K_{12}\times\cdots\times K_{n-1,n},{\rm SO}(n))$ isomorphic to $C(K_{ij}^1,{\rm SO}(3))$, where $K_{ij}^1=K_{ij}\times K_{ik_1}\times K_{k_1j}$. Consequently, $\mathcal{L}_{ij}^1=\overline{{\rm Lie}\{\b_{ij}\O_{ij},\b_{ik_1}\O_{ik_1},\b_{k_1j}\O_{k_1j}\}}$ is isomorphic to $C(K_{ij}^1,\mathfrak{so}(3))$ by Remark \ref{rmk:ensemble_control}. 
Notice that the cardinality of $\{1,\dots,n\}\backslash\{i,j\}$ is $n-2$, so there are $n-2$ distinct subsets of $\mathcal{B}$ (including $\mathcal{S}_1$), denoted by $\mathcal{S}_1,\dots,\mathcal{S}_{n-2}$, in the form of $\mathcal{S}_l=\{\O_{ij},\O_{ik_l},\O_{k_lj}\}$ for some $k_l\in\{1,\dots,n\}\backslash\{i,j\}$, and their intersection only contains $\O_{ij}$. Similar to $\mathcal{L}_{ij}^{1}$, $\mathcal{L}_{ij}^{l}=\overline{{\rm Lie}\{\b_{ij}\O_{ij},\b_{ik_l}\O_{ik_l},\b_{k_lj}\O_{k_lj}\}}$ is isomorphic to $C(K_{ij}^l,\mathfrak{so}(3))$ for each $l=1,\dots, n-2$, where $K_{ij}^l=K_{ij}\times K_{ik_l}\times K_{k_lj}$. As a result, for any $f\in C(K^{\a}_{ij},\mathbb{R})$ and $g\in C(K^{\g}_{ij},\mathbb{R})$ with $\a\neq\g$, we have $f(\b_{ij},\b_{ik_{\a}},\b_{k_{\a}j})\O_{ij}\in\mathcal{L}_{ij}^\a$ and $(g(\b_{ij},\b_{ik_{\g}},\b_{k_{\g}j})/\b_{ik_{\g}})\O_{ik_{\g}}\in\mathcal{L}_{ij}^\b$. Because of 
\begin{align*}
&[[f(\b_{ij},\b_{ik_{\a}},\b_{k_{\a}j})\O_{ij},\b_{ik_{\g}}\O_{ik_{\g}}],(g(\b_{ij},\b_{ik_{\g}},\b_{k_{\g}j})/\b_{ik_{\g}})\O_{ik_{\g}}]\\
&=f(\b_{ij},\b_{ik_{\a}},\b_{k_{\a}j})g(\b_{ij},\b_{ik_{\g}},\b_{k_{\g}j})\O_{ij},
\end{align*}
the set of the coefficients of $\O_{ij}$ in $\overline{{\rm Lie}(\cup_{l=1}^{n-2}\mathcal{L}_{ij}^{l})}$, denoted by $\mathcal{A}_{ij}$, is a subalgebra of $C\big(\prod_{1\leq i<j\leq n}K_{ij},\mathbb{R}\big)$ generated by $C(K_{ij}^1,\mathbb{R}),\dots,C(K_{ij}^{n-2},\mathbb{R})$. Furthermore, let $\mathcal{A}$ denote the subalgebra of $C\big(\prod_{1\leq i<j\leq n}K_{ij},\mathbb{R}\big)$ generated by $\mathcal{A}_{ij}$, $1\leq i<j\leq n$, then $\overline{{\rm Lie}(\cup_{1\leq i<j\leq n}\cup_{l=1}^{n-2}\mathcal{L}_{ij}^{l})}=\mathcal{A}\otimes\mathfrak{so}(n)$ holds. Because $C(K_{ij}^l,\mathbb{R})$ separates points in $K_{ij}^l$ for each $l=1,\dots,n-2$ and $1\leq i<j\leq n$ as shown in the proof of Proposition \ref{prop:so3_3dispersion}, $\mathcal{A}$ is able to separate points in $\prod_{1\leq i<j\leq n}K_{ij}$. By Stone-Weierstrass theorem, $\mathcal{A}$ is dense in $C\big(\prod_{1\leq i<j\leq n}K_{ij},\mathbb{R}\big)$, and then so is $\mathcal{A}\otimes\mathfrak{so}(n)$ in $C\big(\prod_{1\leq i<j\leq n}K_{ij},\mathbb{R}\big)\otimes\mathfrak{so}(n)=C\big(\prod_{1\leq i<j\leq n}K_{ij},\mathfrak{so}(n)\big)$. Notice that $\mathcal{A}\otimes\mathfrak{so}(n)\subseteq\overline{{\rm Lie}(\{\b_{ij}\O_{ij}:1\leq i<j\leq n\})}$ holds by the construction of $\mathcal{A}$, thus we conclude $\overline{{\rm Lie}(\{\b_{ij}\O_{ij}:1\leq i<j\leq n\})}=C\big(\prod_{1\leq i<j\leq n}K_{ij},\mathfrak{so}(n)\big)$, which then implies ensemble controllability of the system in \eqref{eq:son_ensemble_full} on $C\big(\prod_{1\leq i<j\leq n}K_{ij},{\rm SO}(n)\big)$.

(Step II): Given the ensemble system in \eqref{eq:son_ensemble}, there is an ensemble system in the form of \eqref{eq:son_ensemble_full} so that these two systems have the same controllable submanifold.


By the condition that each individual system in the ensemble system \eqref{eq:son_ensemble} is controllable on SO$(n)$, any $\O_{ij}\in\mathcal{B}$ can be generated by iterated Lie brackets of the elements in $\mathcal{F}=\{\O_{i_1j_1},\dots,\O_{i_mj_m}\}$. As a result, for each $\O_{ij}\not\in\F$, there exists a positive monomial function $\eta_{ij}:K\rightarrow\mathbb{H}$ such that $\eta_{ij}(\b)\O_{ij}$ can be generated by sucessively Lie bracketing the elements in $\G=\{\b_{i_1j_1}\O_{i_1j_1},\dots,\b_{i_mj_m}\O_{i_mj_m}\}$. Now, consider the following ensemble system,
\begin{align}
\label{eq:son_ensemble_equivalent}
&\frac{d}{dt}X(t,\b)=\Big[\sum_{\O_{ij}\in\F}\b_{ij}u(t)\O_{ij}+\sum_{\O_{ij}\not\in\F}\eta_{ij}(\b)u_{ij}(t)\O_{ij}\Big]X,\nonumber\\
&X(0,\b)=I,
\end{align}
its controllable submanifold has Lie algebra $\overline{{\rm Lie}(\mathcal{G}\cup\mathcal{G}')}$, where $\mathcal{G}'=\{\eta_{ij}(\b)\O_{ij}:\O_{ij}\not\in\F\}$. Because $\eta_{ij}(\b)\O_{ij}\in\Lie(\G)$ for each $i,j=1,\dots,n$, $\Lie(\G)=\Lie(\G\cup\G')$ holds, which also implies $\overline{\Lie(\G)}=\overline{\Lie(\G\cup\G')}$. Since we have shown that $\overline{{\rm Lie}(\G)}$ is the Lie algebra of the controllable submanifold of the system in \eqref{eq:son_ensemble}, the two ensmeble systems \eqref{eq:son_ensemble} and \eqref{eq:son_ensemble_equivalent} have the same controllable submanifold.

(Step III): The system in \eqref{eq:son_ensemble} is ensemble controllable on $C(K,{\rm SO}(n))$.


In step II, we have shown that each $\eta_{ij}(\b)$ is a positive monomial function defined on the compact subset $K$ of $\mathbb{H}^m$, where we define $\eta_{i_kj_k}(\b)=\b_{i_kj_k}$ for $k=1,\dots,m$. Let $\mathcal{R}_{ij}=\eta_{ij}(K)$ be the image of $\eta_{ij}$, then $\mathcal{R}=\prod_{1\leq i<j\leq n}\mathcal{R}_{ij}$ is a compact subset of $\mathbb{H}^{n(n-1)/2}$ by the continuity of each $\eta_{ij}$ and Tychonoff's product theorem \cite{Munkres00}. Then, the conclusion in Step I implies that the following ensemble system parameterized by $\eta=(\eta_{12},\dots,\eta_{n-1,n})\in\mathcal{R}$ 
\begin{align}
\label{eq:son_enesmble_reparametrized}
\frac{d}{dt}X(t,\eta)=\Big[\sum_{1\leq i<j\leq n}\eta_{ij}v_{ij}(t)\O_{ij}\Big]X(t,\eta),\quad X(0,\eta)=I
\end{align}
is ensemble controllable on $C(\mathcal{R},\mathfrak{so}(n))$. 

Now, consider $\eta$ as a function of $\b$ from $K$ to $\mathcal{R}$ given by $(\b_{i_1j_1},\dots,\b_{i_mj_m})\mapsto(\b_{i_1j_1},\dots,\b_{i_mj_m},\dots,\eta_{n,n-1})$, then $\eta$ is smooth and its differential 
$$d\eta=\left[\begin{array}{c} I_m \\ * \end{array}\right],$$
is full rank, where $I_m$ is the $m$-by-$m$ identity matrix. This implies that $\eta$ is a smooth embedding, and hence $\eta(K)$ is a compact $m$-dimensional embedded submanifold of $\mathcal{R}$ \cite{Lee03}. 
By Tietze's Extension Theorem \cite{Munkres00}, for any $f\in C(\eta(K),{\rm SO}(n))$, there exists $g\in C(\mathcal{R},{\rm SO}(n))$ such that $f=g|\eta(K)$, which implies that the map from $C(\mathcal{R},{\rm SO}(n))$ to $ C(\eta(K),{\rm SO}(n))$ given by $g\mapsto g|\eta(K)$ is surjective. Then, by Step II, ensemble controllability of the system in \eqref{eq:son_enesmble_reparametrized} on $C(\mathcal{R},{\rm SO}(n))$ leads to ensemble controllability of the system in \eqref{eq:son_ensemble} on $C(\eta(K),{\rm SO}(n))$.
Moreover, since $\eta$ is a diffeomorphism between $K$ and $\eta(K)$, the function from $C(K,{\rm SO}(n))$ to $C(\eta(K),{\rm SO}(n))$ given by $f\mapsto f\circ\eta^{-1}$ is a Lie group isomorphism, which then concludes ensemble controllability of the system in \eqref{eq:son_ensemble} on $C(K,{\rm SO}(n))$.\hfill$\Box$

In Step III above, 
the key observation leading to ensemble controllability of the system in \eqref{eq:son_ensemble} is the compactness of $\eta(K)\subset\mathbb{H}^{n(n-1)/2}$. Consequently, the proof still holds if the parameter space is diffeomorphic to a compact submanifold of the upper half space as shown in the following corollary.

\begin{corollary}
\label{cor:son_ensemble_embedding}
The ensemble of systems defined on ${\rm SO}(n)$, given by 
\begin{align}
\label{eq:son_ensemble_embedding}
\frac{d}{dt}X(t,\b)=\Big[\sum_{k=1}^m f_k(\b) u_k(t)\, \Omega_{i_kj_k}\Big]X(t,\b),\quad X(0,\b)=I,
\end{align}
is ensemble controllable on $C(K,{\rm SO}(n))$ if and only if each individual system with respect to a fixed $\b\in K$ in this ensemble is controllable on ${\rm SO}(n)$, where $K$ is a compact smooth manifold, and $f:K\rightarrow\mathbb{H}^m$ defined by $\b\mapsto(f_1(\b),\dots,f_m(\b))$ is a smooth embedding.
\end{corollary}
{\it Proof.}
The necessity is clear, and thus we only need to prove the sufficiency. By defining $\eta_i=f_i(\b)$ for each $i=1,\dots,m$, Theorem \ref{thm:son_ensemble} implies that the system in \eqref{eq:son_ensemble_embedding} parameterized by $\eta=(\eta_1,\dots,\eta_m)'$ is ensemble controllable on $C(f(K),{\rm SO}(n))$. In addition, because $f$ is a smooth embedding, the map from $C(K,{\rm SO}(n))$ to $C(f(K),{\rm SO}(n))$ given by $g\mapsto g\circ f^{-1}$ is a Lie group isomorphism, and hence the system in \eqref{eq:son_ensemble_embedding} is ensemble controllable on $C(K,{\rm SO}(n))$.\hfill$\Box$

Because Step I in the proof of Theorem \ref{thm:son_ensemble} follows from ensemble controllability of systems on SO$(3)$, this theorem, as well as Corollary \ref{cor:son_ensemble_embedding}, exclude systems defined on SO$(2)$. 

\begin{remark}
	\label{rm:so2}
	\rm	
	An ensemble of bilinear systems defined on {\rm SO}$(2)$ is not ensemble controllable. Because $\mathfrak{so}(2)$ is a one-dimensional real vector space with the only basis element $\O_{12}$, any ensemble system on {\rm SO}$(2)$ in the form of \eqref{eq:son_ensemble} can be uniquely represented by
	\begin{align}
		\label{eq:so2}
		\frac{d}{dt}X(t,\b)=\b u(t)\O_{12}X(t,\b)=\b u(t)\left[\begin{array}{cc} 0 & -1 \\ 1 & 0 \end{array}\right]X(t,\b),\quad X(0,\b)=I,
	\end{align}
	where $\beta$ is the parameter taking values on a compact set $K\subset\mathbb{H}$, $X(t,\cdot)\in C(K,{\rm SO}(2))$ is the state, and $u(t)\in\mathbb{R}$ is a piecewise constant control input. However, $\mathfrak{so}(2)$ is {\it nilpotent}, which disables the generation of terms $\b^k\O_{12}$ for $k\geq2$ by iterated Lie brackets of the single control vector field $\b\O_{12}$. As a result, $\overline{{\rm Lie}(\b\O_{12})}$ only contains first order terms 
	of $\b$, and hence the system in \eqref{eq:so2} is ensemble uncontrollable on $C(K,{\rm SO}(2))$.
\end{remark}

\subsection{Ensemble controllability of systems on semisimple Lie groups}
\label{sec:SUn}
The equivalence between ensemble and classical controllability established in Theorem \ref{thm:sen_ensemble} reduced the evaluation of controllability for infinite-dimensional ensemble systems to finite-dimensional single systems. This reduction made it possible to explicitly characterize the generically intractable ensemble controllability property using classical approaches for finite-dimensional control systems, i.e., the LARC for bilinear systems and the Kalman rank condition for linear systems. A natural question concomitant with this property for systems on SO$(n)$ is what other classes of ensemble systems inherit such equivalence in controllability to their subsystems. In this section, we show that ensemble systems defined on semisimple Lie groups exhibit such an equivalence property.

To illuminate this extension, we begin with our discussion on the system defined on SU$(2)$, the special unitary group of $2\times2$ unitary matrices with determinant 1, which is also the most fundamental semisimple Lie group. Notice that its Lie algebra $\mathfrak{su}(2)$, containing all $2\times2$ skew-Hermitian traceless matrices, is isomorphic to $\mathfrak{so}(3)$ by identifying the three basis elements of $\mathfrak{su}(2)$,
\begin{align*}
B_1=\frac{1}{\sqrt{2}}\left[\begin{array}{cc} 0 & i \\ i & 0 \end{array}\right],\quad B_2=\frac{1}{\sqrt{2}}\left[\begin{array}{cc} 0 & -1 \\ 1 & 0 \end{array}\right],\quad\text{and}\quad B_3=\frac{1}{\sqrt{2}}\left[\begin{array}{cc} i & 0 \\ 0 & -i \end{array}\right],
\end{align*}
with $\Omega_x$, $\Omega_y$ and $\Omega_z$, respectively, and $B_1$, $B_2$, and $B_3$ are the Pauli matrices multiplied by $i/\sqrt{2}$, where $i$ is the imaginary unit. In particular, this is called the spin representation of $\mathfrak{su}(2)$. Consequently, following the same proof as that of Proposition \ref{prop:so3_3dispersion}, the system defined on SU(2),
\begin{align*}
\frac{d}{dt}X(t,\beta)=\Big[\sum_{k=1}^3\beta_ku_kB_k\Big]X(t,\beta)
\end{align*} 
is ensemble controllable on $C(K,{\rm SU}(2))$, where $\beta=(\b_1,\b_2,\b_3)$ is the parameter vector taking values on a compact set $K\subset\mathbb{H}^3$. 
This result forms the basis of investigating ensemble controllability for systems evolving on semisimple Lie groups using the covering method. The prerequisite for this investigation is to cover semisimple Lie groups by Lie subgroups isomorphic to SU(2). Similar to Example \ref{ex:cover}, it suffices to construct covers consisting of Lie subalgebras isomorphic to $\mathfrak{su}(2)$.

Given a semisimple Lie group $G$, its semisimple Lie algebra $\mathfrak{g}$ admits a root space decomposition as $\mathfrak{g}=\mathfrak{h}\oplus\bigoplus_{\a\in R}\mathfrak{g}_\a$, where $\mathfrak{h}$ is the Cartan subalgebra, $R$ is the set of nonzero roots, and $\mathfrak{g}_\a$ is the space of root vectors for the root $\a$ \cite{Hall15}. Then, for each root $\a\in\mathbb{R}$, we can construct a Lie subalgebra $\mathfrak{s}_\a$ of $\mathfrak{g}$ so that $\mathfrak{s}_\a$ is isomorphic to $\mathfrak{su}(2)$. To proceed, we first equip the Cartan subalgebra $\mathfrak{h}$ an inner product $\<\cdot,\cdot\>$, through which we define the notion of coroot of $\a$ as $H_\a=2\a/\<\a,\a\>$. Then,  any element $X_\a\in\mathfrak{g}_\a$ satisfies $[H_\a,X_\a]=\<\a,H_\a\>X_\a=2X_\a$ by the definition of a root. Let $Y_\a=-\bar X_\a$, where $\bar X_\a$ denotes the complex conjugate of $X_\a$, then we can show that $Y_\a\in\mathfrak{g}_{-\a}$, $[H_\a,Y_\a]=-2Y_\a$, and $[X_\a,Y_\a]=H_\a$. As a result, $H_\a$, $X_\a$, and $Y_\a$ generate a Lie subalgebra of $\mathfrak{g}$ isomorphic to $\mathfrak{su}(2)$, denoted by $\mathfrak{s}_\a$. However, $H_\a$, $X_\a$ and $Y_\a$ do not give rise to the spin representation of $\mathfrak{s}_\a$ as desired, i.e., $H_\a$, $X_\a$, and $Y_\a$ do not satisfy the same Lie bracket relations as $B_1$, $B_2$ and $B_3$. To construct the spin representation of $\mathfrak{s}_\a$, 
we further define $B_1^\a=i H_\a/2$, $B_2^\a=i(X_\a+Y_\a)/2$ and $B_3^\a=(Y_\a-X_\a)/2$, which lead to the Lie bracket relations $[B_1^\a,B_2^\a]=B_3^\a$, $[B_2^\a,B_3^\a]=B_1^\a$, and $[B_3^\a,B_1^\a]=B_2^\a$. Moreover, because the roots span the Cartan subalgebra $\mathfrak{h}$ \cite{Hall15}, we have constructed a cover of $\mathfrak{g}$ as $\mathcal{U}=\{\mathfrak{s}_\a:\a\in R\}$, in which each $\mathfrak{s}_\a=\Lie(\mathcal{B}^\a)={\rm Lie}(\{B_1^\a,B_2^\a,B_3^\a\})$ is isomorphic to $\mathfrak{su}(2)$ with the spin representation. As a result, the proof of Theorem \ref{thm:son_ensemble} for systems on SO($n$) can be adopted to show ensemble controllability of systems evolving on semisimple Lie groups based on covering its Lie algebra by Lie subalgebras in the form of $\mathfrak{s}_\a$ that are isomorphic to $\mathfrak{su}(2)$ with the spin representation.

\begin{theorem}
	\label{thm:semisimple}
	Given an ensemble of bilinear systems defined on a semisimple Lie group $G$ of the form,
	\begin{align}
	\label{eq:semisimple}
	\frac{d}{dt}X(t,\beta)=\sum_{k=1}^m\Big[\beta_ku_k(t)B_k\Big]X(t,\beta), \quad X(0,\beta)=I,
	\end{align}
	where $\beta=(\beta_1,\dots,\beta_m)$ is the parameter vector taking values on a compact subset $K$ of $\mathbb{H}^m$, $X(t,\cdot)\in C(K,G)$ is the state, $u_k(t)\in\mathbb{R}$ are piecewise constant control inputs, and $I$ denotes the identity element of $G$; $B_1,\ldots,B_m$ are elements in the Lie algebra $\mathfrak{g}$ of $G$ with the property that for any $B_i$, $i=1,\ldots,m$, there exist some $B_j$ and $B_k$ such that the Lie subalgebra of $\mathfrak{g}$ generated by $\{B_i,B_j,B_k\}$ is isomorphic to the spin representation of $\mathfrak{su}(2)$. Then, this system is ensemble controllable on $C(K,G)$ if and only if each individual system with respect to a fixed $\b\in K$ in this ensemble is controllable on $G$.	
\end{theorem}
{\it Proof.} The proof is constructive based on the construction described above and then follow the proof of Theorem \ref{thm:son_ensemble}. To be more specific, after obtaining the cover $\mathcal{U}=\{\mathfrak{s}_\a:\a\in R\}$ of $\mathfrak{g}$, we adopt the proof of Theorem \ref{thm:son_ensemble} by replacing $\mathcal{S}_l=\{\O_{ij},\O_{ik_l},\O_{k_lj}\}$ by $\mathcal{B}^\a=\{B_1^\a,B_2^\a,B_3^\a\}$.
\hfill$\Box$

Note that when the semisimple Lie algebra $\mathfrak{g}$ associated with the system in \eqref{eq:semisimple} is over $\mathbb{C}$, the field of complex numbers, the control inputs $u_k$ are also required to be complexed-valued. Correspondingly, the Lie subalgebra of $\mathfrak{g}$ generated by $\{B_i,B_j,B_k\}$ is the special linear Lie algebra $\mathfrak{sl}(2,\mathbb{C})$, the vector space over $\mathbb{C}$ consisting of 2-by-2 complex matrices with trace 0. This is because $\mathfrak{sl}(2,\mathbb{C})$ is the complexification of $\mathfrak{su}(2)$, that is, for any $A\in\frak{sl}(n,2)$ there exist $A_1,A_2\in\mathfrak{su}(2)$ such that $A=A_1+iA_2$, 
\cite{Hall15}.


\begin{remark}
	\label{rmk:not_semisimple}
	A bilinear ensemble system of the form, 
	$$\frac{d}{dt}X(t,\b)=\Big[\sum_{i=1}^m\b_i \, u_i(t)B_i\Big]X(t,\b),$$
	evolving on a Lie group $G$ that is not semisimple can never be ensemble controllable. To see this, let $\mathfrak{g}$ be the Lie algebra of $G$, then $\mathfrak{g}$ has a nontrivial center $\mathfrak{z}$, whose elements commute with every element in $\mathfrak{g}$. Suppose $B_i\in\mathfrak{z}$ for some $i=1,\dots,m$, then $[\b_iB_i,\b_jB_j]=0$ for any $j=1,\dots,m$. Consequently, the Lie algebra generated by the control vector fields is a module of $\mathfrak{g}$ over a space of functions independent of $\beta_i$, and hence the system cannot be ensemble controllable (on a space of functions of $\b_1$, $\dots$, $\b_m$). 
\end{remark}

\section{Ensemble control of systems defined on non-compact Lie groups}
\label{sec:sen}
In Section \ref{sec:covering}, by introducing the covering method, we established the equivalence between ensemble and classical controllability for parameterized populations of bilinear systems evolving on compact and connected Lie groups. Fortunately, this equivalence also holds true for broader classes of bilinear systems, for example, for bilinear systems induced by Lie group actions on vector spaces. The finding sheds light on possible extension of the equivalence property to systems defined on non-compact Lie groups. In particular, we will show that the system evolving on the special Euclidean group SE$(n)$, which contains the action of SO$(n)$ on $\mathbb{R}^n$, 
inherits this property. Moreover, it is also worth noting that the action of SO$(n)$ on $\mathbb{R}^n$ is neither free nor transitive. In the following section, we briefly review some essential properties of the Lie group SE$(n)$ and its Lie algebra $\mathfrak{se}(n)$ as a prerequisite for carrying out the analysis of ensemble controllability for the system defined on SE$(n)$.

\subsection{Basics of the SE$(n)$ and $\mathfrak{se}(n)$}
Consider the Euclidean space $\mathbb{R}^n$ as a Lie group under addition, then its semidirect product with SO$(n)$, denoted by SE$(n)=\mathbb{R}^n\rtimes {\rm SO}(n)$, is called the special Euclidean group. Therefore, every element in SE$(n)$ can be represented by a 2-tuple $(x,X)$ with $x\in\mathbb{R}^n$ and $X\in{\rm SO}(n)$. Algebraically, the group multiplication is given by $(x,X)(y,Y)=(x+Xy,XY)$ for any $x,y\in\mathbb{R}^n$ and $X,Y\in {\rm SO}(n)$, which also indicates that $(0,I)$ is the identity element of SE$(n)$. Topologically, due to the non-compactness of $\mathbb{R}^n$, SE$(n)$ is also a non-compact Lie group. In addition, SE$(n)$ can be smoothly embedded into GL$(n+1,\mathbb{R})$, the general linear group consisting of all $(n+1)$-by-$(n+1)$ invertible matrices. This embedding immediately yields a matrix representation for each $(x,X)\in {\rm SE}(n)$ as 
$$
(x,X)=\left[\begin{array}{cc} X & x \\ 0 & 1 \end{array}\right],
$$
which also reveals that SE$(n)$ contains SO$(n)$ and $\mathbb{R}^n$ as Lie subgroups. 

Geometrically, let $\gamma(t)=(x(t),X(t))$ be a smooth curve in SE$(n)$ with $\gamma(0)=(0,I)$, then its time derivative at $t=0$, i.e., $\dot{\g}(0)=(\dot{x}(0),\dot{X}(0))$, 
gives rise to an element in the Lie algebra $\mathfrak{se}(n)$ by identifying $\mathfrak{se}(n)$ with $T_{(0,I)}{\rm SE}(n)$, the tangent space of SE$(n)$ at the identity $(0,I)$. Note that $X(t)$ is a curve in SO$(n)$ with $X(0)=I$, and hence we have $\dot{X}(0)\in\mathfrak{so}(n)$. Therefore, every element $(v,\O)\in\mathfrak{se}(n)$  also admits a matrix representation as
\begin{align*}
(v,\O)=\left[\begin{array}{cc} \O & v \\ 0 & 0 \end{array}\right],
\end{align*}
where $\O\in\mathfrak{so}(n)$ and $v\in\mathbb{R}^n$. 

Similar to $\mathfrak{so}(n)$, $\mathfrak{se}(n)$ is also a finite-dimensional vector space, and hence has a basis. Let $\{e_1,\dots,e_n\}$ denote the standard basis of $\mathbb{R}^n$, and define $\mathcal{R}=\{R_{ij}\in\mathfrak{se}(n):R_{ij}=(0,\O_{ij}),1\leq i<j\leq n\}$ and $\mathcal{T}=\{T_{k}\in\mathfrak{se}(n):T_k=(e_k,0),1\leq k\leq n\}$, then the set $\mathcal{R}\cup\mathcal{T}$ forms a basis of $\mathfrak{se}(n)$. The following lemma then characterizes the Lie bracket relations among the basis elements of $\mathfrak{se}(n)$.

\begin{lemma}
\label{lem:sen_lie}
	The Lie brackets among elements in the basis of $\mathfrak{se}(n)$ satisfy that $[R_{ij},R_{kl}]=\d_{jk}R_{il}+\d_{il}R_{jk}+\d_{jl}R_{ki}+\d_{ik}R_{lj}$, $[R_{ij},T_k]=\d_{jk}T_i-\d_{ik}T_j$, and $[T_k,T_l]=0$ for all $1\leq i,j,k,l\leq n$, where $\d$ is the Kronecker delta function.
\end{lemma}
{\it Proof.}
The proof follows from direction computations of Lie brackets by using the matrix representations of $R_{ij}$, $R_{kl}$, $T_k$, and $T_l$.\hfill$\Box$

Notice that Lie brackets among the elements in $\mathcal{R}=\{R_{ij}:1\leq i<j\leq n\}$ follow the same relation as those elements in $\B=\{\Omega_{ij}:1\leq i<j\leq n\}$ as shown in Lemma \ref{lem:son}. This indicates that the Lie algebra $\mathfrak{se}(n)$ contains $\mathfrak{so}(n)$ as a Lie subalgebra. Together with the inclusion of SO$(n)$ in SE$(n)$ as a Lie subgroup, a system defined on SE$(n)$ also contains a system on SO$(n)$ as a subsystem. These relations will help facilitate the controllability analysis of the system on SE$(n)$.


\subsection{A decomposition method for controllability analysis of systems on \bf{\rm SE}$(n)$}
\label{sec:decomposition}
In this section, we focus on the controllability analysis of a single bilinear system defined on SE$(n)$, which builds the foundation towards examining controllability of an ensemble of such systems detailed in the next section. This analysis also illuminates the framework for analyzing controllability of systems induced by Lie group actions on vector spaces. Controllability of systems induced by Lie group actions has been extensively studied \cite{Boothby75,Boothby79,Jurdjevic96}, however, these previous works were largely restricted to consider systems induced by free or transitive Lie group actions. Unfortunately, the action of SE$(n)$ on $\mathbb{R}^n$ is neither free nor transitive, which disables the use of the previously developed conditions to examine controllability of systems on SE$(n)$. Here, we leverage the semidirect product structure of SE$(n)$ to decompose a system defined on this Lie group into two components, the rotational (SO$(n)$) and translational ($\mathbb{R}^n$) components, so that controllability of SE$(n)$ can be analyzed by individually examining that of each component. This approach works for systems on SE$(n)$ because the semidirect product structure is independent of the freeness and transitivity of the group action. It is also potentially applicable to systems induced by general Lie group actions.

For systems on SE$(n)$, we are particularly interested in those governed by the vector fields in $\mathcal{R}\cup\mathcal{T}$ of the form,
\begin{align}
\label{eq:sen}
\frac{d}{dt}\left[\begin{array}{cc} X & x \\ 0 & 1 \end{array}\right]&=\left(\sum_{s=1}^{m_1}u_s(t)\left[\begin{array}{cc}\O_{i_sj_s} & 0 \\ 0 & 0 \end{array}\right]+\sum_{l=1}^{m_2}v_l(t)\left[\begin{array}{cc} 0 & e_{k_l} \\ 0 & 0 \end{array}\right]\right)\left[\begin{array}{cc} X & x \\ 0 & 1 \end{array}\right], \\
(x(0),X(0))&=(0,I),\nonumber
\end{align}
where $\O_{i_sj_s}\in\mathcal{B}$ is a basis element of $\mathfrak{so}(n)$, $e_{k_l}$ is the $k_l$-th standard basis vector of $\mathbb{R}^n$, and $u_s(t),v_l(t)\in\mathbb{R}$ are piecewise constant control functions for all $s=1,\dots,m_1$ and $l=1,\dots,m_2$. 
Because SE$(n)$ contains SO$(n)$ and $\mathbb{R}^n$ as Lie subgroups, the system in \eqref{eq:sen} can be decomposed into two subsystems on SO$(n)$ and $\mathbb{R}^n$, given by
\begin{align}
&\dot X(t)=\Big[\sum_{s=1}^{m_1}u_s(t)\O_{i_sj_s}\Big]X(t),\quad X(0)=I, \label{eq:sen_r}\\
&\dot x(t)=\Big[\sum_{s=1}^{m_1}u_s(t)\O_{i_sj_s}\Big]x(t)+\sum_{l=1}^{m_2}v_l(t)e_{k_l}, \quad x(0)=0, \label{eq:sen_t}
\end{align}
representing the rotational and translational dynamics of the system, respectively. This decomposition enables a tractable way to understand controllability of the system in \eqref{eq:sen}.


\begin{theorem}
\label{thm:sen}
A system defined on {\rm SE}$(n)$ as in \eqref{eq:sen} is controllable if and only if its rotational component in \eqref{eq:sen_r} and translational component in \eqref{eq:sen_t} are simultaneously controllable on {\rm SO}$(n)$ and $\mathbb{R}^n$, respectively.
\end{theorem}
{\it Proof.}
(Necessity): Geometrically, SE$(n)$ is trivially diffeomorphic to $\mathbb{R}^n\times{\rm SO}(n)$ through the identity map $(x,X)\mapsto(x,X)$. Therefore, if the system in \eqref{eq:sen} is controllable on SE$(n)$, then the direct product of the controllable submanifolds of its subsystems in \eqref{eq:sen_t} and \eqref{eq:sen_r} must be $\mathbb{R}^n\times{\rm SO}(n)$, and hence, the systems in \eqref{eq:sen_r} and \eqref{eq:sen_t} are controllable on SO$(n)$ and $\mathbb{R}^n$, respectively. 

(Sufficiency): Given any $X_F\in {\rm SO}(n)$ and $x_F\in\mathbb{R}^n$, it suffices to show that there exist piecewise constant control inputs $u_1,\dots,u_{m_1},v_1,\dots,v_{m_2}$ that simultaneously steer the systems in \eqref{eq:sen_r} from $I$ to $X_F$ and \eqref{eq:sen_t} from $0$ to $x_F$. 

At first, we claim that $m_2\geq1$ must hold if the system in \eqref{eq:sen_t} is controllable on $\mathbb{R}^n$. Otherwise, the system reduces to 
\begin{align}
\label{eq:SOn_on_Rn}
\dot x(t)=\Big[\sum_{s=1}^{m_1}u_s(t)\O_{i_sj_s}\Big]x(t),
\end{align}
which describes the dynamics of the system in \eqref{eq:sen_r} on SO$(n)$ acting on $\mathbb{R}^n$. However, the homogeneous spaces of the Lie group action of SO$(n)$ on $\mathbb{R}^n$ are spheres centered at the origin \citep{Lee03}. Consequently, the controllable submanifold of the system in \eqref{eq:SOn_on_Rn} must be contained in a sphere, which contradicts the controllability of the system on $\mathbb{R}^n$.

Now, let $\mathbb{S}^{n-1}_{\|x_F\|}$ denote the sphere centered at the origin with radius $\|x_F\|$, where $\|\cdot\|$ denotes the Euclidean norm on $\mathbb{R}^n$, and $V$ be the subspace of $\mathbb{R}^n$ spanned by $e_{k_1},\dots,e_{k_{m_2}}$, then $V\cap\mathbb{S}^{n-1}_{\|x_F\|}\neq\varnothing$ holds. Pick a point $z\in V\cap\mathbb{S}^{n-1}_{\|x_F\|}$, because SO$(n)$ acts on $\mathbb{S}^{n-1}_{\|x_F\|}$ transitively \citep{Lee03}, there exists $A\in{\rm SO}(n)$ such that $x_F=Az$.

In the following, we will develop a control strategy to simultaneously steer the system in \eqref{eq:sen_r} from $I$ to $X_F$ and the system in \eqref{eq:sen_t} from $0$ to $x_F$ in three steps. First, because the system in \eqref{eq:sen_r} is controllable on SO$(n)$, the control inputs $u_1,\dots,u_{m_1}$ can be appropriately designed to steer the system  from $I$ to $A^{-1}X_F$, and simultaneously, the system in \eqref{eq:sen_t} stays at the origin by setting $v_1=\cdots=v_{m_2}=0$. Then, we set $u_1=\dots=u_{m_1}=0$ and apply $v_1,\dots,v_{m_2}$ to steer the system in \eqref{eq:sen_t} from the origin to $z$. In this step, the rotational component in \eqref{eq:sen_r} stays at $A^{-1}X_F$. At last,  $u_1,\dots,u_{m_2}$ can be turned on again to steer the system in \eqref{eq:sen_r} from $A^{-1}X_F$ to $X_F$. Since $x_F=Az$, the translational component in \eqref{eq:sen_t} will be simultaneously steered to $x_F$ from $z$, which also completes the proof. \hfill$\Box$

The proof of Theorem \ref{thm:sen} indeed provides a systematic control design procedure to simultaneously steer the systems in \eqref{eq:sen_r} and \eqref{eq:sen_t} between desired states, which concludes controllability of the system in \eqref{eq:sen}. Alternatively, the proof can also be carried out algebraically by computing the Lie algebras generated by the control vector fields of these systems.  
Furthermore, notice that 
the translational component in \eqref{eq:sen_t} also involves the rotational dynamics 
through the SO$(n)$ action on $\mathbb{R}^n$, 
therefore, it is possible to completely determine controllability of the system in \eqref{eq:sen} on SE$(n)$ solely by its translational component in \eqref{eq:sen_t} on $\mathbb{R}^n$.

%
\begin{corollary}
\label{cor:sen}
A system on {\rm SE}$(n)$ as in \eqref{eq:sen} is controllable if and only if its translational component in \eqref{eq:sen_t} is controllable on $\mathbb{R}^n$ and remains controllable on $\mathbb{S}^{n-1}$ if $x(0)\in\mathbb{S}^{n-1}$ and $v_l=0$ for all $l=1,\dots,m_2$, where $\mathbb{S}^{n-1}$ denotes the $(n-1)$-dimensional unit sphere centered at the origin.
\end{corollary}
{\it Proof.} We have shown in the proof of Theorem \ref{thm:sen} that if $v_1=\dots=v_{m_2}=0$, then the rotational component in \eqref{eq:sen_t} reduces to a system induced by the action of SO$(n)$ on $\mathbb{R}^n$. The conclusion then follows from the fact that this Lie group action is transitive on $\mathbb{S}^{n-1}$ \cite{Lee03}. \hfill$\Box$

The above analyses for a single system defined on SE$(n)$ offer the basics for us to move on to the ensemble case in the next section.

\subsection{Ensemble controllability of systems on {\rm SE}$(n)$}
In this section, we will investigate controllability of an ensemble of bilinear systems defined on SE$(n)$. In particular, we focus on the ensemble of the form, 
\begin{align}
	\label{eq:sen_ensemble}
	&\frac{d}{dt}\left[\begin{array}{cc} X(t,\b) & x(t,\b) \\ 0 & 1 \end{array}\right] = \sum_{s=1}^{m_1}u_s(t)\left[\begin{array}{cc} \b_s\O_{i_sj_s} & 0 \\ 0 & 0 \end{array}\right]\left[\begin{array}{cc} X(t,\b) & x(t,\b) \\ 0 & 1 \end{array}\right] \nonumber\\
	&+\sum_{l=1}^{m_2}v_l(t)\left[\begin{array}{cc} 0 & e_{k_l} \\ 0 & 0 \end{array}\right]\left[\begin{array}{cc} X(t,\b) & x(t,\b) \\ 0 & 1 \end{array}\right],\quad X(0,\beta) = I,\quad x(0,\b)=0, 
\end{align}
where $\b=(\b_1,\dots,\b_{m_1})$ is the parameter vector varying on a compact set $K\subset\mathbb{H}^{m_1}$, $\O_{i_sj_s}\in\mathcal{B}$ is a standard basis element of $\mathfrak{so}(n)$ for each $s=1,\dots,m_1$, and $e_{k_l}$ is the $k_l$-th standard basis vector of $\mathbb{R}^n$ for each $l=1,\dots,m_2$. Analogous to the case of a single bilinear system defined on SE$(n)$ discussed in the previous section, the ensemble system in \eqref{eq:sen_ensemble} also admits a decomposition into its rotational and translational components as follows,
\begin{align}
&\frac{d}{dt}X(t,\b)=\Big[\sum_{s=1}^{m_1}\b_su_s(t)\O_{i_sj_s}\Big]X(t,\b),\quad X(0,\b)=I \label{eq:sen_r_ensemble},\\
&\frac{d}{dt}x(t,\b)=\Big[\sum_{s=1}^{m_1}\b_su_s(t)\O_{i_sj_s}\Big]x(t,\b)+\sum_{l=1}^{m_2}v_l(t)e_{k_l},\quad x(0,\b)=0, \label{eq:sen_t_ensemble}
\end{align}
which in turn leads to a characterization of ensemble controllability of the system in \eqref{eq:sen_ensemble} in terms of ensemble controllability of its rotational and translational components in \eqref{eq:sen_r_ensemble} and \eqref{eq:sen_t_ensemble}, respectively.


\begin{theorem}
\label{thm:sen_ensemble}
An ensemble of systems as in \eqref{eq:sen_ensemble} is ensemble controllable on $C(K,{\rm SE}(n))$ if and only if its rotational component in \eqref{eq:sen_r_ensemble} and translational component in \eqref{eq:sen_t_ensemble} are ensemble controllable on $C(K,{\rm SO}(n))$ and $C(K,\mathbb{R}^n)$, respectively.
\end{theorem}
{\it Proof.} The proof is based on the development of a control strategy that simultaneously steers the ensemble systems in \eqref{eq:sen_r_ensemble} and \eqref{eq:sen_t_ensemble} between the respective desired states, which follows the same proof as for Theorem \ref{thm:sen}. Alternatively, we can also adopt the covering method by acting the cover $\mathcal{U}=\{\mathcal{L}_{ij}^l:l=1,\dots,n-2,1\leq i<j\leq n\}$ of $C(K,\mathfrak{so}(n))$ constructed in Theorem \ref{thm:son_ensemble} on $\mathbb{R}^n$. Consequently, $\mathcal{U}\cup\{e_{k_1},\dots,e_{k_{m_2}}\}$ forms a cover of $C(K,\mathbb{R}^n)$, treated as the Lie algebra of the Lie group $C(K,\mathbb{R}^n)$. Then, the rest of the proof follows that of Theorem \ref{thm:son_ensemble}.
\hfill$\Box$


In Theorem \ref{thm:son_ensemble}, we proved the remarkable result that an ensemble system on $C(K,{\rm SO}(n))$ is ensemble controllable if and only if each individual system in this ensemble is controllable on SO$(n)$. By using the decomposition in \eqref{eq:sen_r_ensemble} and \eqref{eq:sen_t_ensemble}, 
this equivalence between ensemble controllability and classical controllability can be extended to ensemble systems defined on $C(K,{\rm SE}(n))$. 

\begin{corollary}
\label{cor:sen_ensemble}
The system in \eqref{eq:sen_ensemble} is ensemble controllable on $C(K,{\rm SE}(n))$ if and only if each individual system in this ensemble is controllable on ${\rm SE}(n)$.
\end{corollary} 
{\it Proof.} To facilitate the proof, we define the notations $\F_1=\{\O_{i_1j_1},\dots,\O_{i_{m_1}j_{m_1}}\}$, $\F_2=\{\O_{i_1j_1}x,\dots,\O_{i_{m_1}j_{m_1}}x,e_{k_1},\dots,e_{k_{m_2}}\}$, $\G_1=\{\b_1\O_{i_1j_1},\dots,\b_{m_1}\O_{i_{m_1}j_{m_1}}\}$, and $\G_2=\{\b_1\O_{i_1j_1}x,\dots,\b_{m_1}\O_{i_{m_1}j_{m_1}}x,e_{k_1},\dots,e_{k_{m_2}}\}$.

The necessity is obvious, so it remains to prove the sufficiency. Assume that each system with a fixed $\b\in K$ in the ensemble \eqref{eq:sen_ensemble} is controllable on SE$(n)$, then by Theorem \ref{thm:sen}, any individual system in the ensemble \eqref{eq:sen_r_ensemble} or \eqref{eq:sen_t_ensemble} is also controllable on SO$(n)$ or $\mathbb{R}^n$, respectively. Hence, the ensemble system in \eqref{eq:sen_r_ensemble} is ensemble controllable on $C(K,{\rm SO}(n))$ by Theorem \ref{thm:son_ensemble}. Then, Theorem \ref{thm:sen_ensemble} implies that it suffices to prove ensemble controllability of the system in \eqref{eq:sen_t_ensemble} on $C(K,\mathbb{R}^n)=C(K,\mathbb{R})\otimes\mathbb{R}^n$, which is equivalent to showing $f(\b)e_k\in\overline{\Lie(\G_2)}$ for any standard basis element $e_k\in\mathbb{R}^n$ and  $f\in C(K,\mathbb{R})$ by Remark \ref{rmk:ensemble_control}.


Because each individual system in the ensemble \eqref{eq:sen_t_ensemble} is controllable on $\mathbb{R}^n$, there exists $\O_{ij}\in\F_1$ and $e_l\in\F_2$ such that $[\O_{ij}x,e_l]=e_k$. 
Furthermore, ensemble controllability of the system in \eqref{eq:sen_r_ensemble} guarantees $f(\b)\O_{ij}\in\overline{\Lie(\G_1)}$, which then gives $[f(\b)\O_{ij}x,e_l]=f(\b)e_k$, i.e., $f(\b)e_k\in\overline{\Lie(\G_2)}$. Therefore, the ensemble system in \eqref{eq:sen_t_ensemble} is ensemble controllable on $C(K,\mathbb{R}^n)$. \hfill$\Box$ 

As a consequence of Theorem \ref{thm:sen_ensemble} and Corollary \ref{cor:sen_ensemble}, the equivalence between ensemble controllability and classical controllability also holds for the translational component of the ensemble system as in \eqref{eq:sen_t_ensemble}. 
This in turn gives rise to a characterization of ensemble controllability of systems on $C(K,{\rm SE}(n))$ solely by their translational components.

\begin{corollary}
\label{cor:sen_ensemble_2}
The system in \eqref{eq:sen_ensemble} is ensemble controllable on $C(K,$ ${\rm SE}(n))$ if and only if its translational component in \eqref{eq:sen_t_ensemble} is ensemble controllable on $C(K,\mathbb{R}^n)$, and remains ensemble controllable on $C(K,\mathbb{S}^{n-1})$ if $x(0,\cdot)\in C(K,\mathbb{S}^{n-1})$ and $v_l=0$ for all $l=1,\dots,m_2$.
\end{corollary}
{\it Proof.} The proof directly follows from Theorem \ref{thm:sen_ensemble} and Corollaries \ref{cor:sen} and \ref{cor:sen_ensemble}.\hfill$\Box$

Notice that the proof of Corollary \ref{cor:sen_ensemble} relies on ensemble controllability of systems evolving on $C(K,{\rm SO}(n))$. Because all the results regarding ensemble controllability of systems on $C(K,{\rm SO}(n))$ established in Section \ref{sec:covering} 
concerned the cases of $n\geq3$, they do not apply to systems defined on $C(K,{\rm SE}(2))$.

\begin{remark} \rm
	\label{rm:se2}
	An ensemble of systems on {\rm SE(2)} in the form of \eqref{eq:sen_ensemble} admits a decomposition,
	\begin{align}
	\frac{d}{dt}X(t,\b)&=\b u(t)\left[\begin{array}{cc} 0 & -1 \\ 1 & 0 \end{array}\right]X(t,\b), \quad\quad\quad\quad\quad\quad X(0,\b)=I, \label{eq:se2_r}\\
	\frac{d}{dt}x(t,\b)&=\b u(t)\left[\begin{array}{cc} 0 & -1 \\ 1 & 0 \end{array}\right]x(t,\b)+\left[\begin{array}{c} 1 \\ 0 \end{array}\right]v(t),\quad x(0,\b)=0,\label{eq:se2_t}
	\end{align}
	where $X(t,\cdot)\in C(K,{\rm SO(2)})$ and $x(t,\cdot)\in C(K,\mathbb{R}^2)$ for each $t\geq0$, and $\b\in K\subset\mathbb{H}$ with $K$ compact. According to Remark \ref{rm:so2}, the rotational component in \eqref{eq:se2_r} is not ensemble controllable on $C(K,{\rm SO}(2))$, or, equivalently, the translational component in \eqref{eq:se2_t} is not ensemble controllable on $C(K,\mathbb{S}^1)$ for 
	$v(t)=0$ and $x(0,\cdot)\in C(K,\mathbb{S}^1)$. This implies 
	uncontrollability of this ensemble on $C(K,{\rm SE}(2))$ by Theorem \ref{thm:sen_ensemble}. However, this 
	does not hinder controllability of the translational component in \eqref{eq:se2_t} on $C(K,\mathbb{R}^2)$. In particular, let $u(t)=1$ be a constant control input, then the ensemble system in \eqref{eq:se2_r} becomes a linear ensemble system with linear parameter variation, studied in our previous work \cite{Li_TAC16}. Because the system matrix $A(\b)=\beta\left[\begin{array}{cc} 0 & -1 \\ 1 & 0 \end{array}\right]$ has disjoint spectra 
i.e., the images of the two eigenvalue functions, $\l_1(\beta)=i\beta$ and $\l_1(\beta)=-i\beta$, are disjoint, this ensemble system representing the translational component is ensemble controllable \cite{Li_TAC16}.


\end{remark}

\begin{remark}
	In our previous work on linear ensemble systems, the equivalence between ensemble controllability and classical controllability requires disjoint spectrum among the system matrices of individual systems \cite{Li_SICON20}. However, for bilinear ensemble systems, the equivalence revealed by utilizing the covering method holds naturally due to their algebraic structure. 
	This finding also indicates that bilinear ensemble systems are easier to be ensemble controllable than linear ensemble systems, which is owing to the nonlinearity in bilinear systems.
\end{remark}


\section{Conclusion} 
In this paper, we propose a unified framework for analyzing ensemble controllability of bilinear ensemble systems defined on semisimple Lie groups. Our main contribution is to develop the covering method that leverages the covering of the state-space Lie group of an ensemble system by its Lie subgroups to enable the controllability analysis of an ensemble through its ensemble subsystems. Exploiting this method, we establish the equivalence between ensemble and classical controllability. This nontrivial property not only reduces the analysis of infinite-dimensional ensemble systems to finite-dimensional single systems, but also empowers the utilization of controllability conditions developed for classical bilinear systems for examining ensemble controllability for bilinear ensemble systems, for example, the LARC and the symmetric group-theoretic controllability conditions in terms of permutation orbits developed in our recent works \cite{Zhang19,Zhang19_IFAC}. Moreover, this equivalence property holds for bilinear ensembles in which the individual systems are defined on non-compact Lie groups, in particular those induced by Lie group actions on vector spaces. This work broadens our understanding of ensemble control systems and opens the door for systematic investigation of fundamental properties of nonlinear ensemble systems. 

\bibliographystyle{siam}
\bibliography{Ensemble_SOn}

\begin{thebibliography}{10}

\bibitem{Apostol74}
{\sc T.~M. Apostol}, {\em Mathematical Analysis}, Addison-Wesley, 2~ed., 1974.

\bibitem{Augier18}
{\sc N.~Augier, U.~Boscain, and M.~Sigalotti}, {\em Adiabatic ensemble control
  of a continuum of quantum systems}, SIAM Journal on Control and Optimization,
  56 (2018), pp.~4045--4068.

\bibitem{Beauchard10}
{\sc K.~Beauchard, J.-M. Coron, and P.~Rouchon}, {\em Controllability issues
  for continuous-spectrum systems and ensemble controllability of bloch
  equations}, Communications in Mathematical Physics, 296 (2010), pp.~525--557.

\bibitem{Becker12}
{\sc A.~Becker and T.~Bretli}, {\em Approximate steering of a unicycle under
  bounded model perturbation using ensemble control}, IEEE Transactions on
  Robotics, 28 (2012), pp.~580--591.

\bibitem{Boothby79}
{\sc W.~Boothby and E.~Wilson}, {\em Determination of the transitivity of
  bilinear systems}, SIAM Journal on Control and Optimization, 17 (1979),
  pp.~212--221.

\bibitem{Boothby75}
{\sc W.~M. Boothby}, {\em A transitivity problem from control theory}, Journal
  of Differential Equations, 17 (1975), pp.~296 -- 307.

\bibitem{Brockett10}
{\sc R.~Brockett}, {\em On the control of a flock by a leader}, Proceedings of
  the Steklov Institute of Mathematics, 268 (2010), pp.~49--57.

\bibitem{Brockett72}
{\sc R.~W. Brockett}, {\em System theory on group manifolds and coset spaces},
  SIAM Journal on Control and Optimization, 10 (1972), pp.~265--284.

\bibitem{Dong14}
{\sc C.~Chen, D.~Dong, R.~Long, I.~R. Petersen, and H.~A. Rabitz}, {\em
  Sampling-based learning control of inhomogeneous quantum ensembles}, Phys.
  Rev. A, 89 (2014), p.~023402.

\bibitem{Chen16_flocks}
{\sc D.~{Chen}, X.~{Liu}, and H.~{Zhang}}, {\em Switching hierarchical
  leadership in coordinated movement of pigeon flocks}, in 2016 35th Chinese
  Control Conference (CCC), 2016, pp.~1158--1163.

\bibitem{XD_Chen19}
{\sc X.~Chen}, {\em Structure theory for ensemble controllability,
  observability, and duality}, Mathematics of Control, Signals, and Systems, 31
  (2019), p.~7.

\bibitem{Chen2020}
\leavevmode\vrule height 2pt depth -1.6pt width 23pt, {\em Ensemble
  observability of bloch equations with unknown population density},
  Automatica, 119 (2020), p.~109057.

\bibitem{Ching2013b}
{\sc S.~Ching and J.~T. Ritt}, {\em Control strategies for underactuated neural
  ensembles driven by optogenetic stimulation.}, Front Neural Circuits, 7
  (2013), p.~54.

\bibitem{Dirr18}
{\sc G.~Dirr and M.~Sch\"{o}nlein}, {\em Uniform and $l^q$-ensemble
  reachability of parameter-dependent linear systems}, 2018.

\bibitem{Dong10}
{\sc D.~Dong and I.~R. Petersen}, {\em Quantum control theory and applications:
  a survey}, IET Control Theory \& Applications, 4 (2010), pp.~2651--2671.

\bibitem{Dong12}
\leavevmode\vrule height 2pt depth -1.6pt width 23pt, {\em Sliding mode control
  of two-level quantum systems}, Automatica, 48 (2012), pp.~725--735.

\bibitem{Folland99}
{\sc G.~B. Folland}, {\em Real Analysis, Modern Techniques and Their
  Applications}, John Wiley \& Sons, Inc., 2nd~ed., 1999.

\bibitem{Glaser98}
{\sc S.~J. Glaser, T.~Schulte-{Herbr\"uggen}, M.~Sieveking, N.~C.~N.
  O.~Schedletzky, O.~W. {S{\o}rensen}, and C.~Griesinger}, {\em Unitary control
  in quantum ensembles, maximizing signal intensity in coherent spectroscopy},
  Science, 280 (1998), pp.~421--424.

\bibitem{Hall15}
{\sc B.~C. Hall}, {\em Lie Groups, Lie Algebras, and Representations}, vol.~222
  of Graduate Texts in Mathematics, Springer International Publishing, 2~ed.,
  2015.

\bibitem{Hatcher02}
{\sc A.~Hatcher}, {\em Algebraic Topology}, Cambridge University Press, New
  York, 2002.

\bibitem{Helmke14}
{\sc U.~Helmke and M.~Schonlein}, {\em Uniform ensemble controllability for
  one-parameter families of time-invariant linear systems}, Systems \& Control
  Letters, 71 (2014), pp.~69--77.

\bibitem{Jurdjevic96}
{\sc V.~Jurdjevic}, {\em Geometric Control Theory}, Cambridge University Press,
  New York, 1996.

\bibitem{Jurdjevic72}
{\sc V.~Jurdjevic and H.~Sussmann}, {\em Control systems on lie groups},
  Journal of differential equations, 12 (1972), pp.~313--329.

\bibitem{Kafashan2015}
{\sc M.~Kafashan and S.~Ching}, {\em Optimal stimulus scheduling for active
  estimation of evoked brain networks}, Journal of Neural Engineering, 12
  (2015), p.~066011.

\bibitem{Khaneja00}
{\sc N.~Khaneja}, {\em Geometric control in classical and quantum systems},
  2000.

\bibitem{Kuritz19}
{\sc K.~{Kuritz}, S.~{Zeng}, and F.~{Allg{\"o}wer}}, {\em Ensemble
  controllability of cellular oscillators}, IEEE Control Systems Letters, 3
  (2019), pp.~296--301.

\bibitem{Lee03}
{\sc J.~M. Lee}, {\em Introduction to Smooth Manifolds}, vol.~218 of Graduate
  Texts in Mathematics, Springer-Verlag New York, 2003.

\bibitem{Li_TAC11}
{\sc J.-S. Li}, {\em Ensemble control of finite-dimensional time-varying linear
  systems}, IEEE Transactions on Automatic Control, 56 (2011), pp.~345--357.

\bibitem{Li_TAC13}
{\sc J.-S. Li, I.~Dasanayake, and J.~Ruths}, {\em Control and synchronization
  of neuron ensembles}, IEEE Transactions on Automatic Control, 58 (2013),
  pp.~1919--1930.

\bibitem{Li_PRA06}
{\sc J.-S. Li and N.~Khaneja}, {\em Control of inhomogeneous quantum
  ensembles}, Physical Review A, 73 (2006), p.~030302.

\bibitem{Li_TAC09}
\leavevmode\vrule height 2pt depth -1.6pt width 23pt, {\em Ensemble control of
  \text{Bloch} equations}, IEEE Transactions on Automatic Control, 54 (2009),
  pp.~528--536.

\bibitem{Li_TAC16}
{\sc J.-S. Li and J.~Qi}, {\em Ensemble control of time-invariant linear
  systems with linear parameter variation}, IEEE Transactions on Automatic
  Control, 61 (2016), pp.~2808 -- 2820.

\bibitem{Li_NatureComm17}
{\sc J.-S. Li, J.~Ruths, and S.~Glaser}, {\em Exact broadband excitation of
  two-level systems by mapping spins to springs}, Nature Communications, 1
  (2017), p.~446.

\bibitem{Li_PNAS11}
{\sc J.-S. Li, J.~Ruths, T.-Y. Yu, H.~Arthanari, and G.~Wagner}, {\em Optimal
  pulse design in quantum control: A unified computational method}, Proceedings
  of the National Academy of Sciences, 108 (2011), pp.~1879--1884.

\bibitem{Li_SICON20}
{\sc J.-S. Li, W.~Zhang, and L.~Tie}, {\em On separating points for ensemble
  controllability}, SIAM Journal on Control and Optimization,  (accepted).

\bibitem{Munkres00}
{\sc J.~R. Munkres}, {\em Topology}, Prentice Hall, 2~ed., 2000.

\bibitem{Petersen16}
{\sc P.~Petersen}, {\em Riemannian Geometry}, vol.~171 of Graduate Texts in
  Mathematics, Springer International Publishing, 3~ed., 2016.

\bibitem{Gong16}
{\sc C.~Phelps, J.~O. Royset, and Q.~Gong}, {\em Optimal control of uncertain
  systems using sample average approximations}, SIAM Journal on Control and
  Optimization, 54 (2016), pp.~1--29.

\bibitem{Li_SICON17}
{\sc S.~Wang and J.-S. Li}, {\em Fixed-endpoint optimal control of bilinear
  ensemble systems}, SIAM Journal on Control and Optimization, 55 (2017),
  pp.~3039--3065.

\bibitem{Li_Automatica18}
\leavevmode\vrule height 2pt depth -1.6pt width 23pt, {\em Free-endpoint
  optimal control of inhomogeneous bilinear ensemble systems}, Automatica, 95
  (2018), pp.~306--315.

\bibitem{Warner83}
{\sc F.~W. Warner}, {\em Foundations of Differentiable Manifolds and Lie
  Groups}, vol.~94 of Graduate Texts in Mathematics, Springer-Verlag New York,
  1~ed., 1983.

\bibitem{Zeng_16_moment}
{\sc S.~Zeng and F.~Allg{\"o}wer}, {\em A moment-based approach to ensemble
  controllability of linear systems}, Systems \& Control Letters, 98 (2016),
  pp.~49--56.

\bibitem{Zeng16}
{\sc S.~Zeng, H.~Ishii, and F.~Allg{\"o}wer}, {\em Sampled observability and
  state estimation of linear discrete ensembles}, IEEE Transactions on
  Automatic Control, 62 (2016), pp.~2406 -- 2418.

\bibitem{zeng2016tac}
{\sc S.~Zeng, S.~Waldherr, C.~Ebenbauer, and F.~Allg\"{o}wer}, {\em Ensemble
  observability of linear systems}, IEEE Transactions on Automatic Control, 61
  (2016), pp.~1452--1465.

\bibitem{Zhang18}
{\sc W.~Zhang and J.-S. Li}, {\em On controllability of time-varying linear
  population systems with parameters in unbounded sets}, Systems \& Control
  Letters, 118 (2018), pp.~94--100.

\bibitem{Zhang19}
{\sc W.~{Zhang} and J.-S. Li}, {\em Analyzing controllability of bilinear
  systems on symmetric groups: Mapping lie brackets to permutations}, IEEE
  Transactions on Automatic Control,  (2019), pp.~1--1.

\bibitem{Zhang19_IFAC}
{\sc W.~Zhang and J.-S. Li}, {\em A symmetric group method for controllability
  characterization of bilinear systems on the special euclidean group},
  vol.~52, 2019, pp.~412 -- 417.
\newblock 11th IFAC Symposium on Nonlinear Control Systems NOLCOS 2019.

\bibitem{Zlotnik12}
{\sc A.~Zlotnik and J.-S. Li}, {\em Optimal entrainment of neural oscillator
  ensembles}, Journal of Neural Engineering, 9 (2012), p.~046015.

\bibitem{Li_ACC12_SVD}
\leavevmode\vrule height 2pt depth -1.6pt width 23pt, {\em Synthesis of optimal
  ensemble controls for linear systems using the singular value decomposition},
  in 2012 American Control conference, Montreal, June 2012.

\bibitem{Li_NatureComm16}
{\sc A.~Zlotnik, R.~Nagao, I.~Z. Kiss, and J.-S. Li}, {\em Phase-selective
  entrainment of nonlinear oscillator ensembles}, Nature Communications, 7
  (2016), p.~10788.

\end{thebibliography}

\end{document}